\documentclass{amsart}
\usepackage[utf8]{inputenc}

\usepackage{graphicx}
\usepackage{subfigure}

\usepackage[centertags]{amsmath}
\usepackage{amssymb}
\usepackage{amsthm}
\usepackage{hyperref}
\hypersetup{
    colorlinks=true,
    linkcolor=blue,
    filecolor=magenta,      
    urlcolor=cyan,
    pdftitle={Overleaf Example},
    pdfpagemode=FullScreen,
    }
    
    \usepackage{soul}
    \usepackage{xcolor}

    \newcommand{\commentOUT}[1]{}

\newcommand{\mbbS}{\mathbb{S}}

\newcommand{\Real}{\mathbb R}
\newcommand{\Realp}{{\mathbb R}^{\scriptscriptstyle{+}}}
\newcommand{\RRealp}{{\mathbb R}^{2\scriptscriptstyle{+}}}

\newcommand{\Natural}{\mathbb N}
\newcommand{\Rational}{\mathbb Q}
\newcommand{\Integer}{\mathbb Z}
\newcommand{\Complex}{\mathbb C}

\newcommand{\mbtt}{\mathbf{t}}

\newcommand{\teta}{\boldsymbol{\ta}}

\newcommand{\Teta}{\boldsymbol{\Theta}}
\newcommand{\hTeta}{\widehat{\Teta}}
\newcommand{\xx}{\boldsymbol{x}}

\newcommand{\bbS}{\mathcal{S}}

\newcommand{\supp}{{\rm supp\,}}

\newcommand{\pr}{\mbox{\sf P}}
\newcommand{\kr}{\mbox{\sf K}}
\newcommand{\tr}{\mbox{\sf T}}

\newcommand{\ex}{{\bf\sf E\,}}               
\newcommand{\dx}{{\bf\sf {D^2}}}               
\newcommand{\kx}{{\bf\sf {Q^2}}}               
\newcommand{\yx}{{\bf\sf {Y}}}
\newcommand{\ux}{{\bf\sf {U}}}

\newcommand{\calb}{{\mathcal B}}

\newcommand{\support}{{{\rm supp}\,}}

\newcommand{\al}{\alpha}                

\newcommand{\ta}{\theta}                

\newcommand{\ra}{\rightarrow}           
\newcommand{\imp}{\Rightarrow}           

\newcommand\opT{\mathbf{T}}
\newcommand\opTt{\mathbf{T_{\scriptscriptstyle\ta}}} 
\newcommand\opTs{\mathbf{T_{\scriptscriptstyle\dag}}}
\newcommand\opTst{\mathbf{T_{\scriptscriptstyle\dag\ta}}}
\newcommand\opTss{\mathbf{T_{\scriptscriptstyle\ddag}}}
\newcommand\opTsst{\mathbf{T_{\scriptscriptstyle\ddag\ta}}}
\newcommand\sta{{\scriptscriptstyle\ta}}
\newcommand\Rab{R_{\scriptscriptstyle{a,b}}}

\newcommand\Rr{R_{\scriptscriptstyle{r}}}
\newcommand\Ror{R_{\scriptscriptstyle{1-r}}}

\newcommand\hpit{\hat{\pi}_{\sta}}

\newtheorem{thm}{Theorem}[section]
\newtheorem{lemma}{Lemma}[section]
\newtheorem{prop}{Proposition}[section]
\newtheorem{coro}{Corollary}[section]
\newtheorem{defn}{Definition}
\newtheorem{exm}{Example}[section]
\newtheorem{exm*}{Example}
\newtheorem{rem}{Remark}[section]

\newtheorem{asm}{Assumption}

\newcommand{\bbGa}{\mathbf \Gamma}

\newcommand{\cM}{\mathcal{M}}

\newcommand{\dblscr}[2]{{{\scriptscriptstyle{#1}}\atop{\scriptscriptstyle{#2}}}}

\newcommand{\spl}{{\scriptscriptstyle{+}}}
\newcommand{\smi}{{\scriptscriptstyle{-}}}
\newcommand{\spm}{{\scriptscriptstyle{\pm}}}

\newcommand{\xpm}{x^{\spm}}

\newcommand{\ypm}{y^{\spm}}

\newcommand{\zpl}{z^{\spl}}
\newcommand{\zmi}{z^{\smi}}

\usepackage[figuresright]{rotating}

\begin{document}

\title[Iterative Random Functions]{On the iterations of some random functions with Lipschitz number one}
\author{Yingdong Lu \and Tomasz Nowicki}
\address{IBM T.J. Watson Research Center, Yorktown Heights, New York, U.S.A.}
\date{}
\begin{abstract}
For the iterations of $\Real\ni x\mapsto |x-\ta|$  random functions with Lipschitz number one, we represent the dynamics as a Markov chain and  prove its convergence under mild conditions. We also de\-mon\-stra\-te that the Wasserstein metric of any two measures will not increase after the corres\-ponding induced ite\-ra\-tions for measures and identify conditions under which a polynomial convergence rate can be achieved in this metric. We also consider  an associated nonlinear operator on the space of probability measures and identify its fixed points through an detailed analysis of their characteristic functions. 
\end{abstract}

\smallskip

\maketitle


\section{Introduction}\label{sec:intro}
Given a Borel probability $\mu$ on $\Realp$ that we shall call \emph{reference measure},  denote by $\Teta\subset \Realp=\{r\in\Real: r\ge 0\}$ its support.
For $\ta\in\Teta$, define a family of functions and their average:  
\begin{align}  
    \opTt:\Realp \to \Realp, &\quad x\mapsto\opTt(x)=|x-\ta|,
 \label{eqn:def opTt}
   \\
       \opT:\Realp \to \Realp, &\quad x\mapsto\opT(x)=\int_{\ta\in\Teta}\opTt(x)\,\mu(d\ta)\,.
 \label{eqn:def opT} 
\end{align}
Let $\ta_k\in\Teta, k=1,2,\ldots$ be an independent and identically $\mu$-distributed sequence.
Then $x_0=x$ and $x_{k+1}= \opT_{\ta_{k}}(x_k)$ define a process 
in  $\Realp$ with (random)  trajectories $\xx=(x_0,x_1,x_2,\dots)$. $x_k, k=0,1,\ldots,$ naturally forms a Markov chain. Let $\pr(x,A)$ denote its transition probability, that is:
\begin{align*}
\pr(x,A)=&\pr_{\Teta} (x,A)=\mu\{\ta\in\Teta:\opTt(x)\in A\}=\mu((x-A)\cup(x+A)).
\end{align*}

The main result in~\cite{doi:10.1137/TPRBAU000047000002000190000001} states that, if $\mu$ has a bounded but not a lattice support, then there exists a unique invariant distribution for the process $\xx$.
This is listed in~\cite{iosifescu09} as an example of iterated function systems with Lipschitz number one, and its convergence under more general condition for $\mu$, as well as convergence rate, are known to be open. 

In this paper, under some mild assumptions, we present the following:
\begin{itemize}
\item 
In Section~\ref{sec:Trajectories}, Theorem~\ref{thm:tajectories}
states that, for a given initial point, depend\-ing on the arithmetic properties of the set of $\ta$'s, the trajectories under the set of maps $\opTt$ either are dense in some intervals or lie on a discrete lattice. 
\label{res-item: Sec 2}
\item In Section~\ref{sec:Markov Chain}, we treat the Markov chain process generated by $\opTt$ \emph{via} $\opTst$, discuss its Harris recurrence and state in Theorem~\ref{thm: inv prob}
that, under mild regularity assumptions on the reference probability measure $\mu$,  
the Markov chain has a unique invariant probability. 
\label{res-item: gen inv prob}
\item
In Section~\ref{sec:Operators}, we shall extend the notion of the maps $\opTt$ and $\opT$ to the operators $\opTst, \opTs$ on $\cM(\Realp)$, the set of measures defined on $\Realp$, see~\eqref{eqn:def pi of teta} and~\eqref{eqn:def avg pi of teta}, and the operators $\opTsst,\opTss$, see~\eqref{eqn:def gammatheta} and~\eqref{eqndef:gammahat}, on the joint measures $\cM(\RRealp)$, where $\RRealp=\Realp\times\Realp$. In addition, we state some basic properties of these operators.
\label{res-item: operators}
\item 
Section~\ref{sec:Wasserstein} contains two results which uses $\opTss$:
For any two probability distributions $\pi$ and $\rho$, Theorem~\ref{thm:Wasserstein decreases} shows that their \emph{Wasserstein distance} decreases under the iterations of the process; we also identify a parameterized family of probability distributions which include ma\-ny commonly used ones including the exponential family and pro\-ba\-bi\-li\-ty distributions with a bounded support, state in Theorem~\ref{thm:Wasserstein poly conv} that the convergence in Wasserstein metric of any initial pro\-ba\-bi\-li\-ty measures from this family to the invariant measure is at least polynomial, and that polynomial rate can be estimated with the pa\-ra\-me\-ters of the probability distribution family.
\label{res-item: Wass decr}
\item
In Section~\ref{sec:inv_prob}, we investigate the family of probability distributions that are invariant under the transformation, with itself as the \emph{reference measure}. We conclude that the probability generating function of such  probability distributions must be in a special form, a linear transform away from a finite Blaschke product. In the case that the Blaschke product is the M\"obius transform, we are able to identify them completely and establish their relationship with the exponential/geo\-met\-ric distributions.
\label{res-item: refer meas is inv}
\end{itemize}

\subsection*{A few notes}
\begin{enumerate}
\item
The \emph{support} of a Borel  measure $\nu$, denoted  $\supp\nu$, 
consists of all points such that any of their open neighborhoods has a positive measure. In other words, it is a complement of the union of all zero measure open sets.  
\item
A function $h$ is \emph{lower semi continuous (l.s.c.)} if for every $x$, we have  $\liminf_{y\to x} h(y) \ge h(x)$, that is there are no jumps down. We shall say that \emph{the measure is l.s.c.} if the functions  $\mu(x+A)$ and $\mu(x-A)$ are l.s.c. for all points $x$ and Borel sets $A$.
\item
We shall call any set  $L=z+w\cdot \Integer$ for some $z,w\in\Realp$ a \emph{(discrete) lattice}. We say that a measure is supported on a lattice if its support is a subset of some lattice.
 \item
     To \emph{exclude trivial cases}, we shall assume that $\Teta$ is not a singleton, as then there is a single  map $\opTt$, and the system is deterministic with iterations converging to the period two map $x\mapsto \ta-x$ on $[0,\ta]$ with~$\opTt^2$ becoming identity. In particular, when the singleton $\ta=0$ then $\opT_0$ itself is identity, $\opT_0(x)\equiv x$, for $x\in\Realp$.
\end{enumerate}

\section{Trajectories of $\opT$}\label{sec:Trajectories}
We shall extend the map $\opTt(x)=|x-\ta|$ to the map of the sets of $\ta\in T\subset\Realp$:  
for $x\in\Realp$, $\opT_T(x)=\{|x-t|: t\in T\}$, and for subsets of $x\in X\subset\Realp$: $\opT_T(X)=\bigcup_{x\in X}\opT_T(x)$. When $\mbtt\in T^{\Integer}$ is a sequence $(t_0,t_1\dots)$, we understand the trajectory 
$\xx=(x_0,x_1,\dots)\in(\Realp)^{\Integer}$ to be the sequence 
\begin{align*}
x_0=\opT^0_{\mbtt}(x)=x, \quad \hbox{and} \quad x_{n+1}=\opT_{\mbtt}^{n+1}
(x)=\opT_{t_n}\circ\opT_{\mbtt}^n(x)=\opT_{t_n}(x_n).
\end{align*}
For $X, T\subset \Realp$, define $\bbS_T^n$ as the set of points visited at least once by some trajectory up to iteration $n$. Therefore, $\bbS^0_TX=X$, $\bbS^{n+1}_TX=\{|y-t|:t\in T, y\in\bbS^n_TX \}$. Let  $\bbS_T X=\bigcup_{n\ge 0}\bbS^{n}_TX$ and  $\mbbS_T X=\bigcap_{k\ge 0}\overline{\bigcup_{n\ge k}\bbS^{n}_T X}$ be the set of accumulation points of the trajectory, they are such points any neighborhood of which was visited infinitely often.
When there is no ambiguity about the set $T$ we shall simplify and write $\bbS(x)=\bbS_{\Teta}(\{x\})$ and $\mbbS(x)=\mbbS_{\Teta}(\{x\})$.

We shall work with a particular set  $\Teta\subset\Realp$, for which we set
\[
\vartheta=\sup\Teta\quad\text{and}\quad \hTeta
=[0,\vartheta)\,,
\]
up to the endpoint the convex hull of $\{0\}\cup\Teta$. It is clear that the results can be scaled, and as $\Teta$ is not trivially a singleton $\{0\}$, we can assume that $1\in\Teta$.

\begin{thm}\label{thm:tajectories}
Let $\Teta\subset \Realp$.
The following conditions are equivalent:
\begin{enumerate}
    \item \label{thm:traj Teta Lattice}
    $\Teta$ is a subset of some discrete lattice $L$
    \item \label{thm:traj bbS not dense}
    for any $x$ the set $\bbS(x)$ is not dense in $\hTeta$,
    \item \label{thm:traj mbbS on Lattice}
    the set $\mbbS(x)\subset (x+L)\cup (-x+L)$. 
\end{enumerate}
    If any of the above condition does not occur then $\hTeta\subset\mbbS(x)$ for any $x\in\Realp$.
\end{thm}
Before proving the Theorem, we introduce some notions and present a techni\-cal proposition.

\subsection*{Dense set and $\epsilon$--covering}
We say that a set $A$ \emph{is $\epsilon$--covering a set $B$} if, for any $b\in B$, there exists an $a\in A$  with $|b-a|<\epsilon$.
A set $A$ is dense in the set $B$ iff it is its  $\epsilon$--cover for every $\epsilon>0$. We say that a map $R$ is $\epsilon$--covering a set $B$ if there exists a trajectory $\{R^n(x)\}$ which $\epsilon$--covers $B$.

\subsection*{Rotations and their trajectories}
We expand the results on the tra\-jec\-tories from~\cite{doi:10.1137/TPRBAU000047000002000190000001}, in particular their Lemma 3.1. We shall use the properties of the rotation maps $\Rr:x\mapsto x+r \mod 1$ on $\Real$, with rotation number $r\in[0,1]$. The name reflects the observation that after identifying the endpoints $\Rr$ represents a \emph{rotation of the circle} by the angle $2\pi r$. For general reference on rotations, please see e.g. Chapter~I in~\cite{de2012one}. The map $\Ror:x\mapsto x-r \mod 1$ is a rotation by the opposite angle $-2\pi r$. Given a rotation of the circle, two trajectories are equal up to a shift by the difference of the initial points.
For any $0<r<1$ the maps $\Rr$ and $\Ror$ provide a $\min(r,1-r)$--covering of $[0,1]$, but in fact more is true.
For an irrational $r\not\in\Rational$, the trajectories are dense in $[0,1]$ (this dates back to Kronecker) and hence for any $\epsilon>0$, $\Rr$ is $\epsilon$--covering $[0,1]$.  
For rational $0<r=\frac{p}{q}<1$ in its simplest form $(p,q)=1$, the trajectories are periodic with period $q$ and they $\frac{1}{q}$--cover $[0,1]$. In particular when $r=\frac{1}{q}$ then $\Rr$ and $\Ror$ are both $\frac{1}{q}$--covering $[0,1]$.

In order to link the properties of the rotations with the trajectories of $\opT$ we shall analyze a subset of $\bbS(x)$ using only two points $0<a<b\in \Teta$. 
For $k<\lfloor\frac{x}{b}\rfloor$, $\opT_b^k(x)=x-kb$ so that eventually the trajectory lands in $[0,b)$.
If $a<x<b$, then  $x\mapsto x-a = b([\frac{x}{b}+1-\frac{a}{b}]\mod{1})=bR_{r}(\frac{x}{b})$, with $r=1-\frac{a}{b}\in(0,1)$ in one step $\opT_a$.
If  $x<a<b$, then  $x\mapsto |x-a|=a-x\mapsto |b-(a-x)|=x+b-a=b(\frac{x}{b}+1-\frac{a}{b})=bR_{r}(\frac{x}{b})$ in two steps $\opT_b\circ\opT_a$. 
By induction $\bbS(x)$ contains the trajectories of (rescaled) rotation
$R_r$,
\begin{align}
\label{eqn:set_relation}
\bbS(x)\supset \bigcup_{\dblscr{a,b\in\Teta}{0<a<b}}\bigcup_{n\ge 0} bR^n_{1-\frac{a}{b}}\left(\frac{\opT^{\lfloor\frac{x}{b}\rfloor}_b(x)}{b}\right).
\end{align}
\commentOUT{The symbol $b R_{\frac{a}{b}}\left(\frac{x}{b}\right)$
is not a trajectory, but just one image. Hence the union is too small. Even further $\bbS(x)$ 
maybe larger than the union of rotation trajectories. It needs an argument that the other inclusion occurs}.

\begin{prop}
\label{prop:when dense}\ 
    \begin{enumerate}
    \item\label{prop: when dense extension}
    If  $\bbS(x)$ is  dense in (\emph{resp. $\epsilon$--covers}) the interval $[0,1)$, \emph{(where $1\in\Teta$ by convention)} then it is dense in 
    (\emph{resp. $\epsilon$--covers}) the interval $\hTeta$. 
    \item\label{prop: when dense bbS to mbbS} 
    If for any $x$, its trajectories $\bbS(x)$ are dense in $\hTeta$ then $\hTeta\subset\overline{\mbbS(x)}=\mbbS(x)$.
    \item\label{prop: when dense zero acc} 
    If $\Teta$ accumulates at $0$, then for any $\ta\in\Teta$, $\bbS(x)$ is dense in $[0,\ta)$.
    \item\label{prop: when dense rational}
    If $a<b\in\Teta$ and $\frac{a}{b}=\frac{p}{q}\in \Rational$ in its simplest form, then $\bbS$ is $\frac{b}{q}$-covering the interval $[0,b]$.
    \item\label{prop: when dense non zero acc}
    If $\Teta$ has an accumulation point $0<z<\infty$, then $\bbS(x)$ is dense in $[0,z)$. For any $\epsilon>0$, $\bbS(x)$ also  $\epsilon$--covers the interval $[0,z)\cup[0,a)$ for some $a>0$ which also contains an element from $\Teta$.
    \item\label{prop: when dense irr}
    If there are $a<b\in\Teta$ with 
    $\frac{a}{b}\not\in\Rational$, then $\bbS(x)$ is dense in $[0,b)$.
\end{enumerate} 
\end{prop}

\begin{proof}
\begin{itemize}
\item[\eqref{prop: when dense extension}] 
    If $\ta=\vartheta$, there is nothing to prove. Otherwise, there exists an $\eta\in(\ta,\vartheta]\cap\Teta$. The image $\opT_\eta(\bbS(x)\cap[0,\ta))$ is dense in (\emph{resp. $\epsilon$--covers}) $(\eta-\ta,\eta)$.  As long as $\ta<\eta-k\ta$, using the iterations of $\opTt$, we obtain the intervals $(\eta-(k+1)\ta, \eta-k\ta)$, where $\bbS(x)$ is dense (\emph{resp. $\epsilon$--covers}) and their union together with $[0,\ta)$ covers $[0,\eta)$. As it works for any $\Teta\ni\eta\le\vartheta =\sup(\Teta)$, the claim is proven.
\item[\eqref{prop: when dense bbS to mbbS}] 
    If $\bbS(x)$ is dense in $[0,\vartheta)$ for any $x$, it is dense also for any $\opT_{\teta}^{n}(x)$, that is the set $\bbS^n(x)$ is dense in $[0,\vartheta)$. After closure, it contains $[0,\vartheta]$ and so does the limit $\mbbS(x)$.
\item[\eqref{prop: when dense zero acc}] 
    If $\Teta$ accumulates at $0$, then there exists a sequence $\Teta\ni\ta_n\searrow 0$, and $\ta_n<\ta_0=\ta$ for $n\ge 1$. The set $\ta_0 R_{\ta_n/\ta_0}(x/\ta_0)$ $a_n$-covers $[0,\ta)$. By \eqref{eqn:set_relation}, $\bbS(x)$ is dense in $[0,\ta)$.
\item[\eqref{prop: when dense rational}]
    It follows from \eqref{eqn:set_relation} as both $R_{\frac{p}{q}}$ and $R_{1-\frac{p}{q}}$ are $\frac{1}{q}$-covering the interval $[0,1]$.
\item[\eqref{prop: when dense non zero acc}] 
        Let $\Teta\ni \ta_n\to z$ monotonically, for any $\epsilon>0$, there exist $m,n$ such that $1-\frac{\epsilon}{2} < \frac{\ta_m}{\ta_n}<1$, and $\ta_n<2z$.
        By rescaling, $\ta_n R_{\ta_m/\ta_n}$ $\epsilon$-covers $[0,\ta_n]$.
        By taking the limit $\ta_n\to z$, $\bbS(x)$ is $\epsilon$--covering $[0,z)$.
        \commentOUT{\emph{The awkward statement with $[0,z)\sup[0,a)$ is added to deal withe the situation when $z$ is the accumulation of decreasing sequences only.}}

\item[\eqref{prop: when dense irr}]
    When $a<b\in \Teta$ with $a/b\not\in \Rational$, then, 
    by the properties of irrational rotations, $\Rab$ has dense orbits in $[0,b]$. 
\end{itemize}
{\hfill$\Box$}\end{proof}

\begin{proof}[of Theorem~\ref{thm:tajectories}]
We know that $\opT_{\teta}^{n}(x)=\pm x+\sum_{i=0}^{n-1}\pm\ta_i$, where the signs depend on the inequalities $\opT_{\teta}^{i}(x)<\ta_i$. \\
\eqref{thm:traj Teta Lattice}$\imp$\eqref{thm:traj bbS not dense}\\
If $\Teta\subset L$, then also $\sum_{i=0}^{n}\pm \ta_i\in L$, which proves a stronger statement that $\bbS(x)\subset (x+L)\cup (-x+L)$, it is not only discrete but also has a special lattice-like shape. \\
\eqref{thm:traj Teta Lattice}$\imp$\eqref{thm:traj mbbS on Lattice}\\
We use the stronger statement from \eqref{thm:traj Teta Lattice}$\imp$\eqref{thm:traj bbS not dense} above.
Because $\bbS(x)$ turns out to be a discrete set, included in $\subset (x+L)\cup (-x+L)$, its accumulation points $\mbbS(x)$
must be a subset of $\bbS(x)$.\\
\eqref{thm:traj mbbS on Lattice}$\imp$\eqref{thm:traj bbS not dense}\\
When $\mbbS(x)$ is discrete so must be $\bbS(x)$, hence $\bbS(x)$ is not dense in an interval.\\
\eqref{thm:traj bbS not dense}$\imp$\eqref{thm:traj Teta Lattice}\\
Assume that $\bbS(x)$ is not dense in the interval $\hTeta$. Then there exists an $\epsilon>0$, such that for any $a<b\in\Teta$ the map $\Rab$ does not $\epsilon$--cover $[0,b)$, otherwise by Proposition~\ref{prop:when dense}.\eqref{prop: when dense extension}
the set $\bbS(x)$ would $\epsilon$--cover $\hTeta$ for any $\epsilon>0$, which would imply density.

By the cases enumerated in Proposition~\ref{prop:when dense}, we know that $\Teta$ must be discrete and co-rational, that is for any $\eta<\ta\in\Teta$, we have $\frac{\eta}{\ta}\in\Rational$.
By (uniform) rescaling, we can assume that $\Teta\subset \Rational$.
If $\Teta$ is finite and rational, then it is contained in a discrete lattice. Otherwise, discrete and infinite $\Teta$ extends to infinity.
Since $q<b/\epsilon$, and $b=\frac{aq}{p}$, so $q<\frac{aq}{p\epsilon}$ or $p<\frac{a}{\epsilon}$. 
Fix 
$a=\min(\Teta\setminus\{0\})$, it exists as $\Teta$ is discrete. 
With fixed $\epsilon$,  $p$ can achieve only a finite number of values, independent of $b$. Let $P$ be a product of all such integers $p$ and the denominator $q_a$ of $a=\frac{p_a}{q_a}$. Then for every $b\in\Teta$, $b\cdot P\in\Integer$, that is $\Teta$ is a subset of a discrete lattice.

Finally, if \eqref{thm:traj bbS not dense} does not occur (and thus also none of \eqref{thm:traj Teta Lattice} and \eqref{thm:traj mbbS on Lattice}) then $\bbS(x)$ is dense in $[0,\vartheta)$ for any $x$. this yields $\mbbS(x)=\overline{\mbbS(x)}$ containing $[0,\vartheta]\supset\hTeta$.

{\hfill$\Box$}\end{proof}

\begin{rem}\label{rem: bound on large steps}
    The interval $\hTeta$ is invariant under any $\opT_{\teta}$. If $\vartheta\in\Teta$ then because $|x-\ta|\le \max(x,\ta)$, also $[0,\vartheta]$ is invariant under any $\opT_{\teta}$.
    If $\opT_{\teta}^N(x)>\vartheta$ then 
    $x>\vartheta$ and $\#\{i<N:\ta_i>\delta\}<\frac{x-\vartheta}{\delta}$.
    As long as $x_i=\opT_{\teta}^i>\vartheta$ then 
    $x_{i+1}=|x_i-\ta_i|=x_i-\ta_i$. When $\#\{i<N:\ta_i>\delta\}=m$ then $\vartheta<\opT_{\teta}^N(x)=x_N<x-m\delta$.

\end{rem}
\begin{rem}\label{rem: intervals}
    If $\Teta$ contains an interval then $\bbS(x)$ contains the interval $\hTeta$.
\end{rem}

\section{Markov chain induced by $\opT$}\label{sec:Markov Chain}

Recall that $\Teta\subset \Realp$ is the support of the \emph{reference measure} $\mu$.
We assume that:
\begin{asm} 
\label{asm:Borel}
The reference measure $\mu$ is Borel on $\Realp$, hence also regular and tight (see e.g.~\cite{halmos1976measure}).
\end{asm}
Sometimes we shall assume additionally that:
\begin{asm}
\label{asm:AC}
 The reference measure $\mu$ is absolutely continuous with respect to the Lebesgue measure and contains an interval in its support.    
\end{asm}
\begin{rem}
\label{rem:property_mu}
Measure $\mu$ of every Borel set can be approximate by the measures of open sets from above and by the measures of compact sets from below.
Under the Assumption~\ref{asm:AC}, the  measure $\mu$ is continuous with respect to translations: for any measurable set $A$, $\lim_{x\to y}\mu(x\pm A)=\mu(y \pm A)$. For the general properties of $\mu$ mentioned above, we refer to~\cite{halmos1976measure}.
\end{rem}

\subsection*{Irreducibility of the Markov chain}
Let $L_T(x,A):=\sum_{n>0}\pr_T^n(x,A)$. 
\begin{lemma}
\label{lem:Npositive}
If the measure $\mu$ is (upper-) regular (open sets condition)  then for any non-empty open $A\subset\hTeta$ we have 
$L_{\Teta}(x,A)=\sum_{N=0}^\infty \pr^N(x,A)>0$.  
\end{lemma}
\begin{proof}
For any open set $A\subset \hTeta=[0,\vartheta)$ ($\vartheta=\sup\Teta$),  
we have $\bbS(x)\cap A\not=\emptyset$, that is there are $\ta_i\in\Theta$, $i=1,\dots,N$ such that for any $\teta$ with initial segment $(\ta_1,\dots,\ta_N,\dots)$, we have $\opT^N_{\teta}(x)\in A$. Taking the preimages, we shall find an open neighborhood $U$ of the cylinder consisting of $\teta$'s with fixed $(\ta_1,\dots,\ta_N)$ such that $\opT_U^N(x)\subset A$. Inside $U$, we find a rectangle  that is an open set $U'=U_1\times\dots U_N\times(\Realp)^{\infty}$, 
where $U_i$ are open neighborhoods of $\ta_i$. 
Therefore 
\[
\pr^N(x,A)=\mu^N(\teta:\opT^N_{\teta}(x)\in A)\ge \prod_{i=1}^N \mu(U_i)\cdot 1>0\,.
\]
{\vskip-0.25cm\hfill$\Box$}\end{proof}

\begin{defn}[Continuous component]
If $a$ is a sampling distribution (i.e. a distribution on nonnegative integers) and there exists a substochastic tran\-si\-tion kernel $\tr$ satisfying, 
\begin{align*}
\kr_a(x,A) \ge \tr(x,A), \quad x\in X,\ A\in \calb(X),
\end{align*}
where $\kr_a(x, A):= \sum_{n=0}^\infty P^n(x,A) a(n)$, and $\tr(\cdot, A)$ is lower semicontinuous function for any $A\in\calb(X)$, then $\tr$ is called a continuous component of $\kr_a$. 
\end{defn}
\begin{defn}[$\tr$-chain]
If $\Psi$  is a Markov chain for which there exists a sampling distribution $a$ such that $\kr_a$ possesses a continuous component $\tr$, with $\tr(x,X)>0$ for all $x$, then $\Phi$ is called a $\tr$-chain. 
\end{defn}
From the proof of Lemma~\ref{lem:Npositive}, we have,
\begin{lemma}
\label{lem:T-chain}
The Markov chain $x_n$ is a $\tr$-chain.
\end{lemma}

Referring to Theorem 6.0.1. from Meyn and Tweedie \cite{meyn93} we arrived to:
\begin{coro}
    The Markov chain generated by $\mu$ and $\opT$ is an \emph{irreducible  $\tr$-chain}. Therefore every compact set is \emph{petite}.
\end{coro}


\subsection*{Positive Harris Recurrence of the Markov Chain}

For fixed $y\in \Real^+$, denote:
\begin{align*}
m(y)&=\mu[0,y)=\int_{[0,y)}\,\mu(d\ta), \quad 0\le m\le 1, \\
    \yx(y)&=\ex_{\mu}(y)=\int_{[0,y)}\ta\,\mu(d\ta), \quad
    \ex=\ex_\mu=\ex_\mu(\infty), \quad 0\le \yx\le\ex, \\
\kx&=\kx_\mu=\int_{[0,\infty)}\ta^2\,\mu(d\ta), \quad
    \dx=\dx_\mu=\kx_\mu-(\ex_\mu)^2.    
\end{align*}
\begin{lemma}\label{lem:Ym}
    For any distribution $\mu$ on $\Realp$ we have for every $x,y\in\Realp$
    
    \begin{align*}
        x\le \frac{\yx(y)-\yx(x)}{m(y)-m(x)}\le y\,,
    \end{align*}
and in particular $
        \frac{\yx(y)}{m}\le y\le \frac{\ex-\yx(y)}{1-m}\,.
$
It follows that $\lim_{y\to\infty}y(1-m(y))\to 0$.
\end{lemma}
\begin{proof}
Using the definitions of $\yx(y)$ and $m$,  we see that
$\yx(y)-\yx(x)=\int_x^y t\mu(dt)\le y\int_x^y\mu(dt)=y\,(m(y)-m(x))$. Similarly $\yx(y)-\yx(x)\ge x\,(m(y)-m(x))$.  So
$\ex-\yx(y)=\int_y^\infty  t\mu(dt)\ge y\int_y^\infty\mu(dt)=y(1-m(y))$.
The expression under the limit is non-negative and bounded from above, since $\ex-\yx(y)\to 0$ as $y\to\infty$.
{\hfill$\Box$}\end{proof}

\begin{figure}
    \centering
    \subfigure[]{\includegraphics[width=0.32\textwidth]{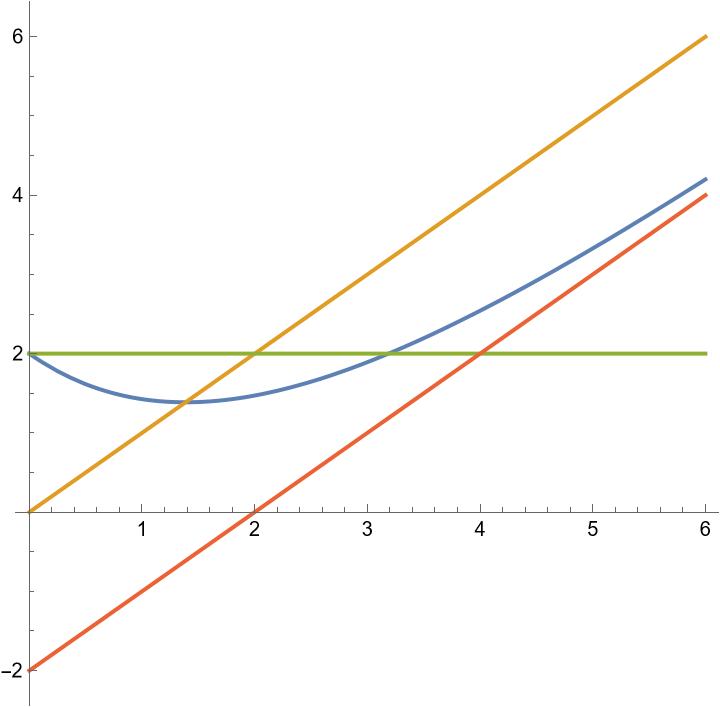}} 
    \subfigure[]{\includegraphics[width=0.32\textwidth]{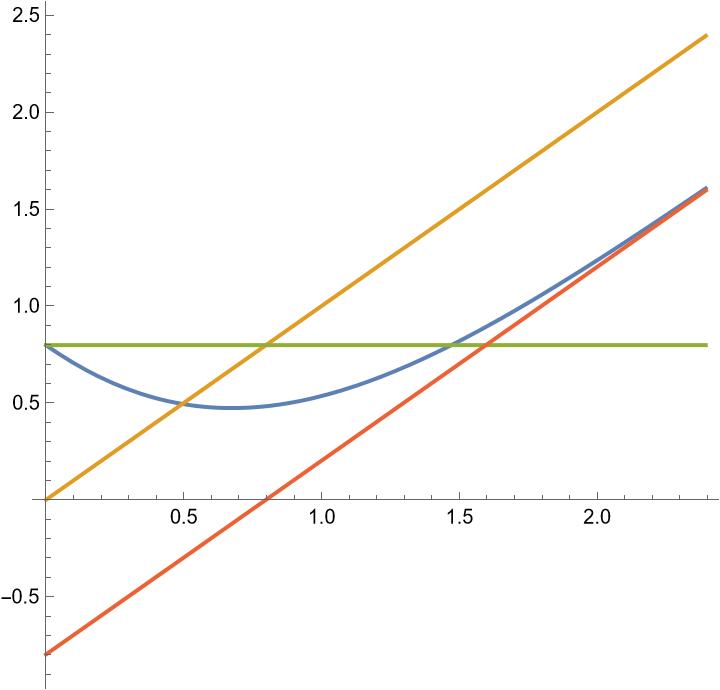}} 
    \subfigure[]{\includegraphics[width=0.32\textwidth]{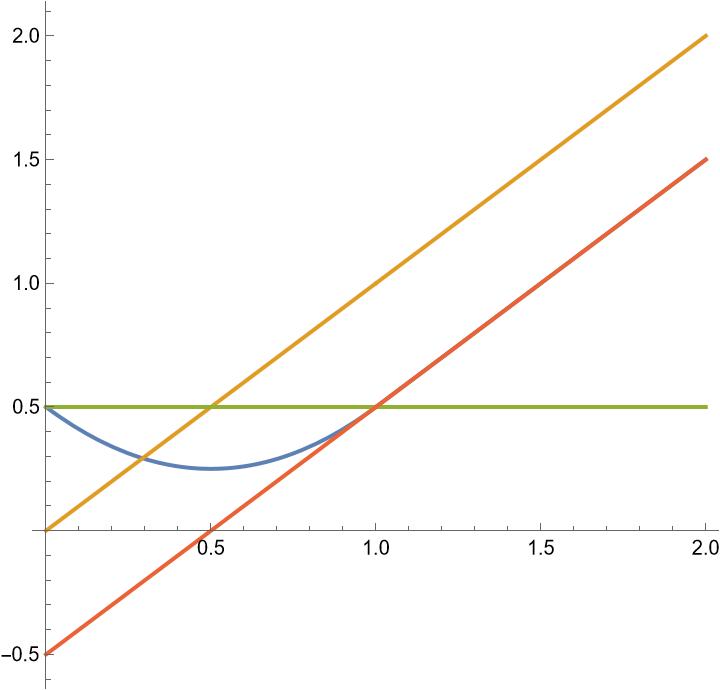}}
    \caption{The $\ux(y)$ function for different reference distributions: (a) exponential (b) normal (one-sided) and  (c) uniform. The horizontal line is at the level $\ex_\mu$, two diagonal lines are identity and the asymptote $y-\ex_\mu$. The fixed point $\ux(y)=y$ lies to the left of $\ex_\mu$, after that $y-\ex_\mu\le \ux(y)<y$.}
    \label{fig:Ufunction}
\end{figure}

Define $\ux(y):=\ex_\mu \opTt(y)$. 
\begin{lemma}\label{lem:property of V}
    The function
    \[
    \ux(y)=y(2m(y)-1)+\ex-2\yx(y)\,,
    \]
    is continuous, 
    $y-\ux(y)\ra \ex$ as $y\to\infty$. 
    Its slope at $y=0$ is  equal to $2\mu(\{0\})-1$. It has a minimum at the median of $\mu$ and it is convex.
\end{lemma}
See Figure~\ref{fig:Ufunction} for $\ux(y)$ function with different reference distributions. 
\begin{proof}
\begin{align*}
\ux(y)&=\ex_\mu[\opTt(y)]= \int_0^y (y-\ta) \mu(\ta) \,d\ta + \int_y^\infty (\ta-y) \mu(\ta) \,d\ta 
\\&= 
y \left[\int_0^y  \mu(\ta) \,d\ta- \int_y^\infty   \mu(\ta) \,d\ta\right]+\left[\int_y^\infty \ta\mu(\ta) \,d\ta -\int_0^y \ta \mu(\ta) \,d\ta \right]
\\&= 
y \left[2\int_0^y\mu(\ta)-\int_0^\infty  \mu(\ta) \,d\ta  \,d\ta\right]+\left[\int_0^\infty\ta \mu(\ta) \,d\ta  - 2\int_0^y \ta \mu(\ta) \,d\ta \right]
\\&=
y\left[2m-1\right]+\left[\ex-2\yx(y)\right] =\left[ym-\yx(y)\right]+\left[(\ex-\yx(y))-(1-m)y\right]\,.
\end{align*}
Clearly $\ux(0)=\ex$.
As 
\[
0\le\lim_{x\searrow 0} \frac{Y(x)}{x}= \lim_{x\searrow 0}\int_{[0,x)}\frac{t}{x}\mu(dt)\le \lim_{x\searrow 0}\frac{0\mu(\{0\})}{x}+\lim_{x\searrow 0}\mu((0,x))=0
\]
we have $\lim_{x\searrow 0}\frac{\ux(x)-\ex}{x}=-1+2\mu(\{0\})$.

Continuity follows from convexity. To prove it directly it is enough to prove the continuity of $I(x)=xm(x)-\yx(x)=x\int_0^x \mu(dt)-\int_0^x t\mu(dt)$.
$I(x+h)-I(x)=h\int_0^{x+h}\mu(dt)+\int_x^{x+h}(x-t)\mu(dt)$. 
Both terms are positive and bounded by $h$ so that $I(x+h)-I(x)\to 0$ as $h\to 0$. 

By Lemma~\ref{lem:Ym}, $0\le \ux(y)-(y-\ex)=2(y(m-1)+\ex-\ex)\to 0$ as $y\to\infty$. 

Let $y$ be such that $m(y)=\frac{1}{2}$. If $x<y$ then 
$\ux(x)-\ux(y)=2xm(x)-x+\ex-2\yx(x)-\ex+2\yx(y)\ge 2xm(x)-x+ 2x (m(y)-m(x))=0$.
Similarly, if $z>y$ then $\ux(z)-\ux(y)=2zm(z)-z-2 \yx(z)+2\yx(y)\ge 
2 zm(z)-z - 2z(m(z)-m(y))=0$.

For $0\le x<y<z$, define 
$$D(x,y):=\frac{\ux(y)-\ux(x)}{y-x}$$ 
and $$R(x,y,z):=\frac
{(y-x)(z-y)}{2}(D(y,z)-D(x,y)),$$ 
then 
\begin{align*}
R(x,y,z)=&   
(y- x) [z m(z)-y m(y)-\yx(z)+ \yx(y)]-\\ 
&- 
(z-y) [y m(y)-x m(x) - \yx(y) + \yx(x)]
\\
=&
(y - x) [z (m(z) - m(y)) -(\yx(z) - \yx(y))]+\\ 
&
+ 
(z-y) [(\yx(y) - \yx(x))- x (m(y) - m(x))].
\end{align*}
Again by Lemma~\ref{lem:Ym}, as $x<y<z$, we have $\yx(y)-\yx(x)\ge x(m(y)-m(x))$ and $\yx(z)-\yx(y)\le z(m(z)-m(y)) $ thus $R(x,y,z)\ge 0$ and the slopes $D$ are increasing for adjacent intervals, which proves convexity of $U$.

The proof is simpler 
when the measure has a density $\mu(dt)=\mu(t)dt$.
Then  $m(y)'=\mu(y)$ and $\yx(y)'=y\mu(y)$
we have $U'(y)=2m(y)-1+2y\mu(y)-2y\mu(y)=2m(y)-1$ and $U''(y)=2\mu(y)\ge 0$, hence convexity and the minimum at the median. Moreover $\lim_{y\to\infty}U'(y)=2-1=1$ and $\lim_{y\to\infty}(\ux(y)-y)=\lim_{y\to\infty}2y(m(y)-1)+\ex-2\ex=-\ex$ by Lemma~\ref{lem:Ym}, hence the asymptote. 
We remark also that $\ux(0)=\ex$ and $U'(0)=-1$.
{\hfill$\Box$}\end{proof}

For $\alpha\in[-1,1)$, let $y_\alpha$ be such that $\ux(y_\alpha)=y_\alpha-\alpha\ex$. We have $y_{-1}=0$, as $\ux(0)=0-(-1)\ex$, and $y_0$ is the unique fixed point $\ux(y_0)=y_0$. 
For $-1\le \alpha< \beta<1$, set the affine function joining the points $(y_\alpha, \ux(y_\alpha))$ and $(y_\beta, \ux(y_\beta))$.
For $y_\alpha\le y\le y_\beta$, let
\[
 \ell_{\scriptscriptstyle \alpha,\beta}(y)=(y_\alpha-\alpha\ex)\frac{y_\beta-y}{y_\beta -y_\alpha}+(y_\beta-\beta\ex)\frac{y-y_\alpha}{y_\beta -y_\alpha}\,,
 \]
we thus have $\ell_{\scriptscriptstyle \alpha,\beta}(y)=y-\frac{\alpha(y_\beta-y)+\beta(y-y_\alpha)}{y_\beta-y_\alpha}\ex$.
Using the asymptote, we can extend the definition of $y_\alpha$ to $\alpha=1$, setting $y_1=\infty$ and extend the definition of $\ell$ to
\[
\ell_{\scriptscriptstyle \alpha,1}(y)=y-\alpha \ex\,.
\]
We see that for $\alpha\ge 0$ and 
$y_\alpha\le y\le y_\beta<1$, we have $\ell_{\scriptscriptstyle \alpha,\beta}(y)\le y$.
\begin{prop}\label{prop: ex going down}
For any $-1\le \alpha<1$, the point $y_\alpha$ exists and is unique. The function $\alpha\mapsto y_\alpha$ is continuous and strictly increasing from $0$ to $+\infty$. Moreover, for 
$-1\le \alpha<\beta\le 1$ and $y\in[y_\alpha,y_\beta)$, we have 
\[
\ux(y)\le \ell_{\scriptscriptstyle \alpha,\beta}(y)\,.
\]
In particular, there exists an $\alpha^*>0$ for which  $y^*=y_{\alpha^*}=(1+\alpha^*)\ex$, so that $\ux(y^*)=\ex$ and $\ux(y)\le \ex$ for $y<y^*$, and $\ux(y)\le y-\alpha^*E$ for $y\ge y^*$.
\end{prop}
\begin{proof}
The proof and the formulae follow from the convexity of $\ux$ and the existence of the asymptote. We remark that for $z\in\{ y_\alpha,y_\beta\}$, we have $\ux(z)=\ell_{\scriptscriptstyle \alpha,\beta}(z)$.

    The function $\ux(y)$ cannot cross twice any line of the family $y\mapsto y-\alpha \ex$ for any $|\alpha|< 1$. First, $\ux(y)$ decreases from the value $\ex$, crosses the diagonal (before or after the minimum) and increases, passing the value $\ex$ for the second time, and then approaches the asymptote $y-\ex$ from above. Thus the uniqueness of $y_\alpha$.
    The inequality follows from $y_\alpha<y<y_\beta$.
{\hfill$\Box$}\end{proof}

\begin{coro}\label{coro: expected value}
   Taking $\alpha=-1$ and $\beta=0$, we have $y_{-1}=0$ and $y_0=q$ the fixed point of $\ux$. For $0\le y\le q$,
   \[
   \ux(y)\le y-\frac{-(q-y)}{q}\ex=y+\ex-\frac{y}{q}\ex\,,
   \]
   and for $0\le \alpha\le y\beta$,
   \[
   \ux(y)\le y-\frac{\alpha(y_\beta-y)+\beta(y-y_\alpha)}{y_\beta-y_\alpha}\ex\le y-\alpha \ex.
   \]
\end{coro}

Foster-Lyapunov Theorem, see, e.g.~\cite{meyn_tweedie_survey,meyn_tweedie_1992, meyn93, Down1997}, states that if there is a function $V(y): \Realp\to\Realp$, a compact set $C$ and constant $\epsilon >0, b\ge 0$, such that, 
\begin{align}
\label{eqn:drift}
\ex [V(x_{n+1})|V(x_n)=y]\le y-\epsilon + b{\bf 1}_{C}(y), 
\end{align}
holds for any $y\in \Realp$, then the Markov chain $x_n$ is positive Harris recurrent, thus has a unique invariant measure.

Corollary~\ref{coro: expected value} implies that the function
$\ux(y)$ we discussed above satisfies the condition~\eqref{eqn:drift}. Hence,

\begin{thm}\label{thm: inv prob}
    Under Assumptions~\ref{asm:Borel} and~\ref{asm:AC}, the Markov chain ge\-ne\-ra\-ted by $\mu$ and $\opT$ has a unique invariant probability.
    {\hfill$\Box$}
\end{thm}


\section{Operators on the space of measures}
\label{sec:Operators}
In this section we shall ex\-tend the notion of the operators acting on $\Realp$ defined in Section~\ref{sec:intro}, namely $\opTt$ by \eqref{eqn:def opTt} for $\ta\in\Teta$ and $\opT$ by \eqref{eqn:def opT},  to the operators defined on the sets of finite measures $\opTst,\opTs:\cM(\Realp)\to\cM(\Realp)$ and, with $\bbGa=\cM(\RRealp)$, 
$\opTsst,\opTss:\bbGa\to\bbGa$.
\begin{align}
\opTst:&\pi\mapsto\hpit; 
&\hpit=(\opTt)_{\#}\pi
\label{eqn:def pi of teta}
\\
\opTs:&\pi\mapsto\hat{\pi};
&\hat{\pi}=\ex_\mu\opTst\pi\,.
\label{eqn:def avg pi of teta}
\end{align}
with the pushforward measure $(\opTt)_{\#}\pi$ defined by $(\opTt)_{\#}\pi(A)=\pi(\opTt^{-1}(A))$ for every measurable set on $\Realp$.
In the case that $\mu$ and $\pi$ have densities (with respect to the Lebesgue measure), abusing slightly the notation with $\mu(dx)=\mu(x)\,dx$ and $\pi(dx)=\pi(x)\,dx$, we can write,
\begin{align*}
\opTst\pi(x)\,dx&
=[\pi(\ta+x)+\pi(\ta-x)]\,dx,
\\
\opTs\pi(x)\,dx&=\left[\int_{\Realp}(\pi(x+\ta)+\pi(\ta-x))\,\mu(d\ta)\right]\,dx,
\end{align*}
recall that $\pi(z),\hpit(z), \hat{\pi}(z)=0$ for $z<0$.


\begin{exm}[Exponential distribution]
\label{exm:expon}
Let $\mu(\ta)d\ta=Me^{-M\ta}d\ta$ with  para\-meter $M>0$. Then, for any exponential distribution with parameter $P$, $\pi(x)dx=P e^{-Px}dx$, we have, 
$\opTs_\mu\pi(x)dx=(\frac{M}{M+P}\mu(x)+\frac{P}{M+P}\pi(x))dx$.
A geometric convergence of the ite\-ra\-tions is evident. 
For starting points $\pi$ having the exponential distribution there is a unique limit distribution, the reference probability $\mu$ itself.  

The uniform convergence extends to the finite convex combinations of exponential densities. For infinite convex combinations there is convergence, but one looses the uniformity. 
{\hfill$\blacksquare$}\end{exm}

\begin{exm}[Bernoulli]\label{exm:Bernoulli}
Let $\mu=p_0 \delta_0+p_1\delta_1$ be supported on two points $0$ and $1$. Then (a) the iterates of $\opTs_\mu$ on any $\pi$ decrease to 0 for $x>1$, and (b) for any $\pi$ with the support on $[0,1]$ the iterates converge to $\frac{\pi(x)+\pi(1-x)}{2}$.

For (a) there is a general argument, see Remark~\ref{rem: bound on large steps}. Let the support of $\mu$ be bounded by (a minimal) $M$ with $\mu[\frac{M}{2},M]>0$. For $z>M$, consider the tail measure: $\opT_\mu\pi[z,\infty)=\pi[z,\infty)-\int_0^M \int_z^{z+\ta}\pi(dw)\mu(d\ta)$. The second term is lower bounded by  $\int_{\frac{M}{2}}^M \int_z^{z+\frac{M}{2}}\pi(dw)\mu(d\ta)=
\int_z^{z+\frac{M}{2}}\pi(dw)\mu[\frac{M}{2},M]$. 
Therefore, for any $z>M$, as long as $\pi[z, z+\frac{M}{2}]>0$, the tail measure decreases. 

For (b) $\opTs_\mu\pi(x)=p_0\pi(x)+p_1\pi(1-x)$ and $\opTs_\mu\pi(1-x)=p_0\pi(1-x)+p_1\pi(x)$. For any $x$ this reduces to the multiplication by a $((p_0,p_1), (p_1, p_0))$ matrix which for the eigenvalue 1 has the eigenvector $(\frac{1}{2},\frac{1}{2})$, the limit of the process. If $\pi$ was supported on $[0,1]$, the limit is $\pi_\infty(x)=\frac{1}{2}(\pi(x)+\pi(1-x))$.
Otherwise, for $x\in[0,1]$: $\pi_\infty(x)=\frac{1}{2}
\sum\limits_{y\,{\rm mod}\,1=x}(\pi(y)+\pi(1-y))$.

In this example for any $\pi$ the limit of iterations exists, but it depends on the starting measure $\pi$. 
{\hfill$\blacksquare$}\end{exm}
We shall denote the measures on $\RRealp$ by $\gamma\in\bbGa$.
There we define the operators  $\opTsst:\gamma\mapsto\hat{\gamma}$ and their averages: 
\begin{align}
\opTsst:\gamma\mapsto\hat{\gamma}_\ta;\quad&\quad\hat{\gamma}_\ta=(\opTt\times\opT)_{\#} \gamma,
\label{eqn:def gammatheta}
\\
\opTss:\gamma\mapsto\hat{\gamma};\quad&\quad\hat{\gamma}=\ex_\mu\opTsst.
\label{eqndef:gammahat}
\end{align}

Finally, denote by  $\Gamma(\rho,\pi)\subset\bbGa$ the set of measures with marginals $\rho$ and $\pi$,  that is for $\gamma\in\Gamma(\rho,\pi)$, we have, $\int_{\Realp}\gamma(x,y)\,dx=\pi(y)$ and $\int_{\Realp}\gamma(x,y)\,dy=\rho(x)$.

Given $\ta\in\Teta$ and $z\in\Realp$, we shall denote
\[
\zpl=\ta+z,\quad \zmi=\ta-z\,.
\]
Obviously, $\{\zpl, \zmi\}=\opTt^{-1}(z)$.

It can be easily verified that $\opTst\pi$ and $\opTs\pi$ are probability measures given that $\pi$ is a probability measure on $\Realp$. Similarly, $\hat{\gamma}_{\ta}$ and $\hat{\gamma}$ define two probability measure, given
$\gamma\in\Gamma(\rho,\pi)$.

\begin{prop}\label{prop:gamma hat}
 For any $\ta$ and any real function $f:\Realp\to\Real$, we have  
 \begin{align*}
&\iint_{\RRealp}f(|x-y|){\hat{\gamma}_\ta} (x,y)(dx,dy)
=
\iint_{\RRealp}f(|\opT _{\ta}x -\opT _{\ta}y|)\,\gamma(dx,dy),
 \\
 &\iint_{\RRealp} f(|x- y|)\hat{\gamma}(dx,dy)
=
\int_{\Teta}\left(\iint_{\RRealp}f(|\opTt x-\opTt y|)\,\gamma(dx,dy))\right)\,\mu(d\ta).
 \end{align*}
\end{prop}
\begin{proof}
We calculate separately the double integrals of the four terms of ${\hat{\gamma}_\ta} (x,y)$. Note that as $z=\opTt(\zpl)=\opTt(\zmi)$, for $x,y\in\Realp$, we always have $|x-y|=|\opTt(\xpm)-\opTt(\ypm)|$.
\begin{align*}
     \iint\limits_{\RRealp}\!f(|x-y|)\gamma(\ta+x,\ta+y)(dx,dy)
     &= \iint\limits_{x^{\spl}, y^{\spl}\ge \ta} \!\!f(|x^{\spl}-y^{\spl}|)\,\gamma(dx^{\spl},dy^{\spl})\\
     &=\iint\limits_{x^{\spl},y^{\spl}\ge \ta}\!\!f(|\opTt x^{\spl}-\opTt y^{\spl}|)\,\gamma(dx^{\spl},dy^{\spl}).
\end{align*}
    Similarly,
  \begin{align*}
    \iint\limits_{\RRealp}f(|x-y|)\,\gamma(\ta-x,\ta+y)(dx,dy)
    &=\!\!\!\iint\limits_{x^{\smi}<\ta\le y^{\spl}}
    \!\!\!f(|\opTt x^{\smi}-\opTt y^{\spl}|)\,\gamma(dx^{\smi},dy^{\spl}),
    \\
    \iint\limits_{\RRealp}f(|x-y|)\,\gamma(\ta+x,\ta-y)(dx,dy)
    &=\!\!\!\iint\limits_{x^{\spl}\ge \ta>y^{\smi}}
    \!\!\!f(|\opTt x^{\smi}-\opTt y^{\spl}|)\,\gamma(dx^{\smi},dy^{\spl}),
    \\
     \iint\limits_{\RRealp}f(|x-y|)\,\gamma(\ta-x,\ta-y)(dx,dy)
    &=\!\!\!\iint\limits_{x^{\smi},y^{\smi}<\ta}
    \!\!\!f(|\opTt x^{\smi}-\opTt y^{\smi}|)\,\gamma(dx^{\smi},dy^{\smi})\,.
    \end{align*}  
After renaming the variables with superscripts to the variables without them,  the sum of the four terms is equal to $\iint_{\RRealp}f(|\opTt x -\opTt y|)\,\gamma(dx,dy)$.    
\begin{align*}
\iint\limits_{\RRealp} f(|x- y|)\hat{\gamma}(dx,dy)&=\iint_{\RRealp}f(|x-y|)\left(\int_{\Teta}\hat{\gamma}_\ta(x,y)\,\mu(d\ta)\right)(dx,dy)\\
&=
\int_{\Teta}\left(
\iint_{\RRealp}f(|x-y|)\hat{\gamma}_\ta(x,y)(dx,dy)\right)\,\mu(d\ta)
\\
&=
\int_{\Teta}\left(\iint_{\RRealp}f(|\opTt x-\opTt y|)\,\gamma(dx,dy))\right)\,\mu(d\ta).\\
\end{align*}
{\vskip-0.25cm\hfill$\Box$}\end{proof}

\begin{figure}
    \centering
    \includegraphics[width=0.60\textwidth]{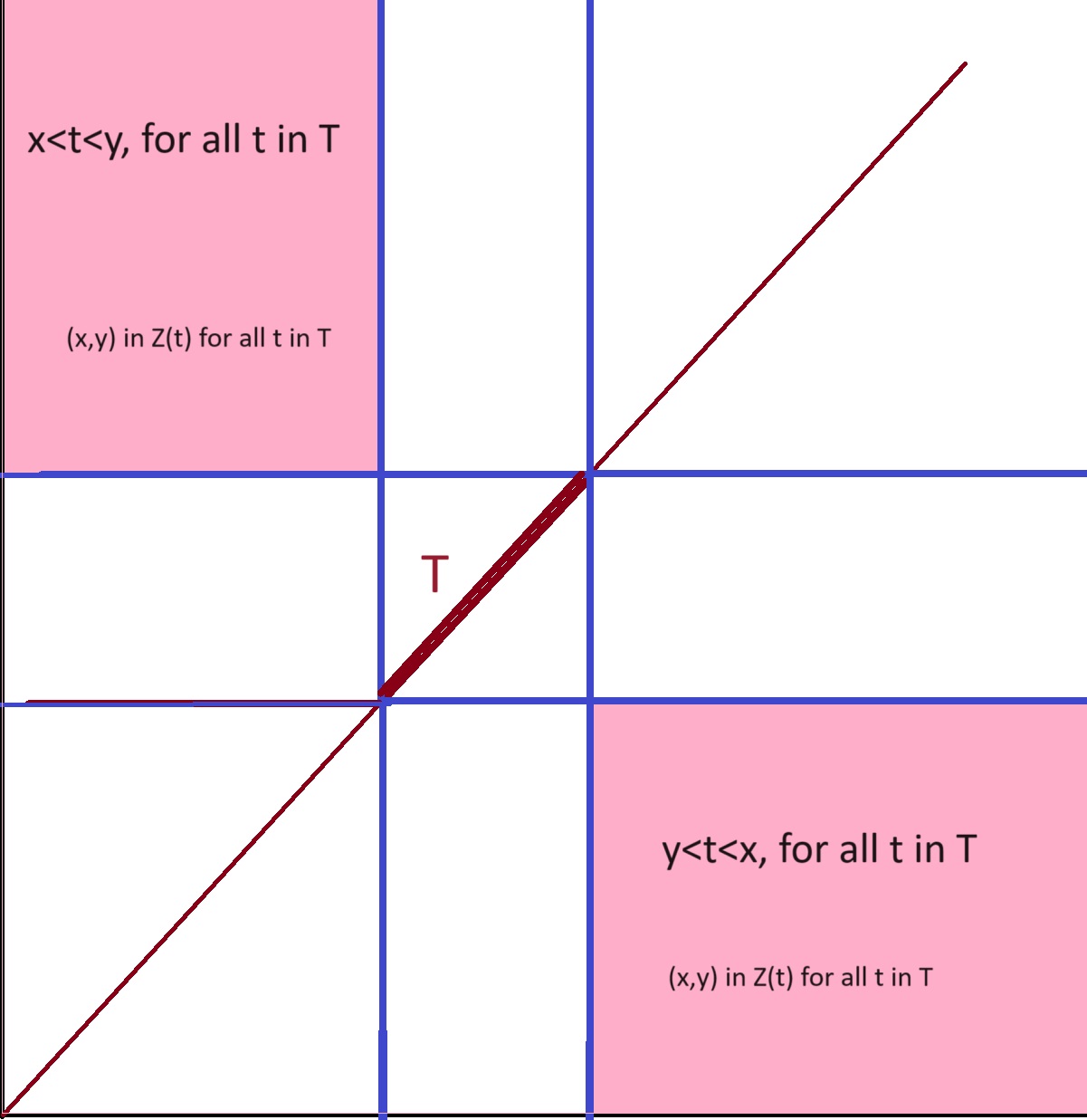}
    \caption{To illustrate $Z(t)$, we put $x$ on the horizontal axis $y$ on the vertical axis, and $t$ on the diagonal $y=x$. $Z(t)$ consists of two domains bounded by the axes and the lines $x=t$ and $y=t$, with interiors disjoint from the diagonal.}
    \label{fig:ZofT}
\end{figure}

Define 
\begin{equation}\label{eqndef:Z of theta}
Z(\ta)=\{(x,y):{\min(x,y)< \ta< \max(x,y)}\}.
\end{equation}
We see that for any $\ta$, the set $Z(\ta)$ is disjoint from the diagonal $x=y$, see Figure~\ref{fig:ZofT}. 

\begin{lemma}\label{lem:subtraction term}
\begin{align*}
     &\iint_{\RRealp} f(|x- y|)\hat{\gamma}(dx,dy)=
     \iint_{\RRealp} f(|x- y|)\gamma(dx,dy)\ - \\
     &-\int_{\ta\in\Teta}\iint_{Z(\ta)}\left(f(|y-x|)-f(|x+y-2\ta|)\right)\gamma(dx,dy)\,\mu(d\ta).
\end{align*}
\end{lemma}
\begin{proof}
    When $\ta\le\min(x,y)$ or $\max(x,y)\le \ta$, $|\opTt(x)-\opTt(y)|=|x-y|$. The remaining cases satisfy 
    $\min(x,y)< \ta < \max(x,y)$, that is $(x,y)\in Z(\ta)$.

{\hfill{$\Box$}}\end{proof}

\begin{coro}\label{coro:gamma hat on non-decr f}
    If $f$ is nondecreasing, then:
    \[
\iint_{\RRealp} f(|x- y|)\,\hat{\gamma}(dx,dy)\le 
\iint_{\RRealp}f(|x-y|)\,\gamma(dx,dy)\,.
    \]
    In particular, it applies to $f(z)=|z|^p$, $p\ge 0$, 
 \[
\iint_{\RRealp}|x-y|^p{\hat{\gamma}}(dx,dy)\le
\iint_{\RRealp}|x-y|^p\,\gamma(dx,dy)\,.
 \] 
\end{coro}
\begin{proof}
We observe that, in $Z(\ta)$, $|x+y-2\ta|=|(x-\ta)-(\ta-y)|\le |x-y|$, thus Lemma~\ref{lem:subtraction term} can be applied. Alternatively, $|\opTt x-\opTt y|\le |x-y|$ implies that $f(|\opTt x-\opTt y|)\le f(|x-y|)$. 
After applying it to Proposition~\ref{prop:gamma hat}, we have: 
\begin{align*}
\iint_{\RRealp} f(|x- y|)\hat{\gamma}(dx,dy)&\le \int_{\Teta}\left(
\iint_{\RRealp}f(|x-y|)\,\gamma(x,y)(dx,dy)\right)\,\mu(d\ta)\\
&=\int_{\Teta}\,\mu(d\ta)\cdot\iint_{\RRealp}f(|x-y|)\,\gamma(dx,dy).
\end{align*}
{\vskip-0.25cm\hfill$\Box$}\end{proof}

\section{The Wasserstein distance $W_p$}\label{sec:Wasserstein}
For a measure $\gamma\in\bbGa$ on $\RRealp=\Realp\times\Realp$, the standard distance ${\rm{dist}}(x,y)=|x-y|$ and $p\ge 1$,  
define 
\[
W_{p}(\gamma):=\left[\iint_{\RRealp} |x-y|^p \gamma(d(x,y))\right]^{\frac{1}{p}}\,.
\]
The \emph{$p$-th Wasserstein distance} on the set measures $\cM(\Realp)$ is given by 
\[
W_p(\rho,\pi)=\inf_{\gamma\in\Gamma(\rho,\pi)}W_{p}(\gamma)\,,
\] 
we recall that $\Gamma(\rho,\pi)\subset\bbGa$ is the set of measures with marginals $\rho$ and $\pi$.
\begin{prop}\label{prop:Wasserstein decreases}
    \begin{align*}
        W_p^p(\opTs\rho,\opTs\pi)\le W_p^p(\rho,\pi).
    \end{align*}
\end{prop}
\begin{proof}
 For any $\epsilon>0$, there exists a $\gamma\in\Gamma(\rho,\pi)$ such that 
 $W_p^p(\gamma)\le W_p^p(\rho,\pi) +\epsilon$. 
Since $\hat{\gamma}=\opTss\gamma\in \Gamma(\opTs\rho, \opTs\pi)$, by Pro\-po\-si\-tion~\ref{prop:gamma hat},  we have
$W_p^p(\hat{\gamma})\le W_p^p(\gamma)$. 
Therefore, 
$$W_p^p(\opTs\rho, \opTs\pi)\le W_p^p(\hat{\gamma})\le W_p^p({\gamma})\le
   W_p^p(\rho,\pi)+\epsilon. $$
{\hfill$\Box$}\end{proof}

\begin{prop}\label{prop:Wgamma strictly decreases}
    If for $Z(\ta)$ defined in~\eqref{eqndef:Z of theta}, there is a $\gamma\in \bbGa$ such that $\int_{\ta\in\Teta}\gamma(Z(\ta))\,\mu(d\ta) >0$,
    then 
    $W_p(\opTss\gamma)<W_p(\gamma)$.
\end{prop}
\begin{rem}    In particular, when $W_p(\rho,\pi)=W_p(\gamma)$ for some $\gamma\in\Gamma(\rho,\pi)$ which satisfies the condition of Proposition~\ref{prop:Wgamma strictly decreases}, $W_p(\opTs\pi, \opTs\rho)<W_p({\gamma})$ strictly.
\end{rem}
For $(x,y)\in Z(\ta)$, set $w(x,y,\ta)=1-\frac{|x+y-2\ta|}{|y-x|}$.
\begin{lemma}
\label{lem:claim}
We have $w=w(x,y,\ta) \in (0,1]$ and $|x+y-2\ta|=|y-x|(1-w)$. 
\end{lemma}
\begin{proof}
$Z(\ta)$ can be partitioned into four subsets depending on the order of the points $x,y$ and the points $\frac{x+y}{2},\ta$ which lie between $x$ and $y$. 
When $x<\ta\le\frac{x+y}{2}<y$ then 
$|x+y-2\ta|=|y-x-2(\ta-x)|=|y-x|(1-2\frac{\ta-x}{y-x})$ where $0<w=2\frac{\ta-x}{y-x}\le 1$,
and similarly for the other three cases.
{\hfill$\Box$}\end{proof}
\begin{proof}[of Proposition~\ref{prop:Wgamma strictly decreases}]
Consider the subtracting term in Lemma~\ref{lem:subtraction term} for $f(z)=z^p$, $p>1$. For such  $(x,y)\in Z(\ta)$, we have $|y-x|^p-|x+y-2\ta|^p=|y-x|^p(1-(1-w)^p)$ for a $w=w(x,y,\ta)\in(0,1]$ from Lemma~\ref{lem:claim}. Meanwhile, we have $1-(1-w)^p\ge w>0$ for $p\ge 1$ by concavity. 
\begin{align*}
&\int_{\ta\in\Teta}\iint_{Z(\ta)}\left(|y-x|^p-|x+y-2\ta|^p\right)\gamma(dx,dy)\,\mu(d\ta)\\
= &\int_{\ta\in\Teta}\iint_{Z(\ta)}|y-x|^p\left(1-(1-w(x,y,\ta))^p\right)\gamma(dx,dy)\,\mu(d\ta)\\
\ge &\int_{\ta\in\Teta}\iint_{Z(\ta)}|y-x|^p w(x,y,\ta)\gamma(dx,dy)\,\mu(d\ta)>0\,,
\end{align*}
as $\gamma(Z(\ta))>0$, $|y-x|>0$ in $Z(\ta)$ and $w>0$.{\hfill$\Box$}\end{proof}

To identify when the inequality in Proposition~\ref{prop:Wasserstein decreases} distance is sharp,  we need the following  technical Lemma:
\begin{lemma}\label{lem:technical A}
   If  $\pi\not=\rho$, $W(\rho,\pi)=W(\gamma)$, then there exist: a set $A\subset \RRealp$ with $\gamma (A)>0$, an $\alpha=\alpha\in(0,1)$  and a proper interval $T=T_\alpha\subset \Realp$, such that $A\subset \bigcap_{s\in T} Z(s)$, and for any $(x,y)\in A$ and any $t\in T$ there holds
   \[
   |x+y-2t|<(1-\alpha)|y-x|\,.
   \]
\end{lemma}
The inclusion condition means that for every $(x,y)\in A$ and every $t\in T$ we have $\min(x,y)\le t\le \max(x,y)$.
\begin{proof}
 As $\pi\not=\rho$, the support of $\gamma$ does not lie entirely on the diagonal $y=x$, in other words there exists an $\epsilon>0$ such that $\gamma(|y-x|\ge \epsilon)>0$. Without the loss of generality, we may assume that $y>x$, so that $\gamma(y-x\ge \epsilon)>0$.
    Given an arbitrary $K>1$, let $\delta=\frac{\epsilon}{K}$. Splitting the set $y-x>\epsilon$ into bands parallel to the diagonal, we see that there exists a $k\ge K$ such that $\gamma(y-x\in[k,k+1)\delta)>0$. Cutting this band by vertical bands, we conclude that there exists $j\ge 0$,  such that the set $$A=\{(x,y): y-x\in[k,k+1)\delta\ \text{and}\ x\in[j,j+1)\delta\}$$ has positive measure $\gamma(A)>0$. 
    In $A$, we have $y\in[k+j,k+j+2)\delta$.
    
    Let $\alpha=\alpha_K\in(0,\frac{K-1}{K+1})\subset(0,1)$, as $K>1$. Then, $k\ge K>\frac{1+\alpha}{1-\alpha}$.
    Define: 
    \begin{align*}
    {t_x}:=&(j+1+\alpha\frac{k+1}{2})\delta\,, {t_y}:=(j+k-\alpha\frac{k+1}{2})\delta\,,  
    T:=T_\alpha=[t_x,t_y].
    \end{align*}
    We have, 
    $t_x>(j+1)\delta\ge x$, $t_y< (j+k)\delta\le y$, $t_y-t_x=k-1-\alpha(k+1)=k(1-\alpha)-(1+\alpha)>0$ and thus $x<t_x<t_y<y$. 

    Suppose first $T\ni t\le \frac{x+y}{2}$. Then, 
    $$|x+y-2t|=(x+y)-2t=(y-x)-2(t-x).$$ We have,
    $${t-x}>{t_x-x}\ge (j+1+\alpha\frac{k+1}{2})\delta-(j+1)\delta=\alpha\frac{k+1}{2}\delta\ge \alpha\frac{y-x}{2}.$$ Thus,
    $|x+y-2t|<|y-x|(1-\alpha)$, since $y-x\ge 0$.
    
    If $T\ni t\ge \frac{x+y}{2}$,
    then $$|x+y-2t|=2t -(y+x)=(y-x)-2(y-t).$$ 
    In this case,
    $${y-t}\ge {y-t_y}\ge (j+k)\delta-(j+k-\alpha\frac{k+1}{2})\delta=\alpha\frac{k+1}{2}\delta> \alpha\frac{y-x}{2}.$$
    As before, $|x+y-2t|<(y-x)-\alpha(y-x)=|y-x|(1-\alpha)$.
    
{\hfill$\Box$}\end{proof}

\begin{thm}\label{thm:Wasserstein decreases}
Suppose the support of $\mu$ is the whole $\Realp$.
If $\pi\not=\rho$ then $W(\opTs\pi,\opTs\rho)<W(\pi,\rho)$ (the inequality is sharp), and the sequence $$W(\opTs^n\pi, \opTs^n\rho)$$ is strictly decreasing with $n$.
\end{thm}
\begin{proof}
    The subtraction term  in Lemma~\ref{lem:subtraction term}  can be estimated from below using Lemma~\ref{lem:technical A}.
    We have $A\subset\RRealp$ with $\gamma(A)>0$, $0<\alpha<1$ and by assumption on the support of $\mu$, $\mu(T)>0$.
    In the proof of Lemma~\ref{lem:technical A}, we have $|y-x|>\epsilon>0$ for $(x,y)\in A$.
    Therefore for strictly increasing functions $f(z)=z^p$, we have:
    \begin{align*}
     &\int_{\Teta}\left(\iint_{Z(\ta)}(f(|y-x|)-f(|x+y-2\ta|)\gamma(dx,dy)\right) \mu(d\ta)
     \\
     &>\int_{\ta\in T}\left(\iint_{A}\left(|y-x|^p-(1-\alpha)^p|y-x|^p\right)\gamma(dx,dy)\right)\mu(d\ta)\\
     \quad&     \ge \epsilon^p(1-(1-\alpha)^p)\gamma(A)\mu(T)>0.
    \end{align*}       
    {\vskip-0.25cm\hfill$\Box$}\end{proof}

\begin{rem}  
    The Wasserstein distance is not decreasing in two cases. 
    \begin{itemize}
    \item
    When $\gamma$ is supported on the vertical and horizontal bands which intersected with the diagonal omit the support of $\mu$, 
    then for almost all $\gamma(x,y)$ both $\mu(x)=0$ and $\mu(y)=0$. The marginals have support disjoint from $\mu$. However by Lemma~\ref{lem:Npositive}, there is an iterate 
    $m$ such that both images of $\opTs^n\rho$ and $\opTs^n\pi$, $n\ge m$, have supports intersecting with the support of $\mu$.
\item
    When all $\gamma$ is supported on the diagonal $x=y$ then  
    the Wasserstein distance is not strictly decreasing. In this case, we have, $$W_p(\opTs\pi, \opTs\rho)=W_p(\pi, \rho),$$ that is the (marginal) projections  of $\gamma$, $\rho$ and $\pi$ have the same distri\-bution.   
    \end{itemize}
\end{rem}

\begin{lemma}
    If the support of $\mu$ contains any interval $[0,2a]\subset \support(\mu)$, then for any $\pi$ there exists an $n$ such that  $\mu(\support((\opTs)^n\pi))>0$. 
    \commentOUT{actually due to Lemma~\ref{lem:Npositive} it is true for any $\mu$ as some open intervals have positive measures.}
\end{lemma}
\begin{proof}
    Consider an interval $[x,y)$ which has positive $\pi$ measure. Then, if $x,y>2a$, using only $\ta\in[a,2a]$ of positive $\mu$ measure one will obtain an interval $[x-2a,y-a]$ of positive $\opTs\pi$ measure. After a few steps we get $y-na<2a$, thus proving the Lemma.  
{\hfill$\Box$}\end{proof}

\begin{coro}
    If the support of $\mu$ contains an initial interval $[0,2a]$,  then the Wasserstein distance decreases strictly after a finite number of steps of $\opTs$.
\end{coro}

\subsection*{Polynomial convergence in Wasserstein metric}
\begin{defn}
A measure $\mu$ on $\Real_+$ is said to \emph{satisfy the $(A,B,C)$-condition} for real numbers $A<1$ and $B, C>0$ if for any interval $(x, y)\subset \Real_+$ such that  $\mu(x,y)>0$, there exists and interval $(L, U)\subset (x,y)$, with $L(x,y)=L< U= U(x,y)$ such that 
\begin{itemize}
\item
for any $z\in [L(x,y),U(x,y)]$ there holds  
$|2z-x-y|\le A|y-x|$, and 
\item 
$\mu ((L,U)) \ge C\mu((x,y))^B$
.\end{itemize}
\end{defn}

The $(A,B,C)$-condition is not very restrictive. In the examples below, we can see that many commonly used distributions satisfy the condition, and identify the constants. 
\begin{exm}
\label{exm:uniform}
Uniform distribution $\mu={\rm Unif}(0,a)$ for some $a>0$ satisfy the con\-di\-tion $(A,B,C)$ with $A= 1-2\kappa$, $B=1$ and $C=1-2\kappa$ for some $\kappa\in (0,\frac12)$; 
It is easy to verify that we can have $L(x,y) = x+\kappa(y-x)$ and $U(x,y) = y-\kappa(y-x)$ . 
{\hfill$\blacksquare$}\end{exm}
\begin{exm}
\label{exm:bounded} Suppose that $\mu$ is a distribution with con\-ti\-nu\-ous density $f$ and a bounded support $[a,b]$. If $f$ is bounded and  $0<m\le f \le M<\infty$, then again, we can verify that 
$A= 1-2\kappa$, $B=1$ and $C=1-2\kappa)m/M$ with $L(x,y) = x+\kappa(y-x)$ and $U(x,y) = y-\kappa(y-x)$.
{\hfill$\blacksquare$}\end{exm}

\begin{exm}
\label{exm:exponential}
Let us first consider the case of exponential distribution with rate one. Let $c>0$, first consider the case $y-x\ge 2c$. Set $L(x,y) = x+c$ for some $c>0$ and $U(x,y) = y-\kappa(y-x)$ for some $\kappa\in (0,\frac12)$ satisfies $(1-\kappa) (x-y)+c <-1$. We have, 
\begin{align*}
\frac{\mu(L,U)}{\mu(x,y)} =& \frac{e^{-(x+c)}- e^{-y+\kappa (y-x)}}{e^{-x}-e^{-y}}\\
=& \frac{e^{-c}- e^{(x-y)-\kappa (x-y)}}{1-e^{x-y}}
\\
=& e^{-c}\frac{1- e^{(1-\kappa) (x-y)+c}}{1-e^{x-y}}.
\end{align*}
The condition $(1-\kappa) (x-y)+c <-1$ implies that $\frac{\mu(L,U)}{\mu(x,y)}$ will be lower bounded by $e^{-c}(1-e^{-1})$. In the case of $y-x< 2c$, let $L(x,y) = x+\kappa(y-x)$ and $U(x,y) = y-\kappa(y-x)$
\begin{align*}
\frac{\mu(L,U)}{\mu(x,y)} =& \frac{e^{-(x+\kappa(y-x))}- e^{-y+\kappa (y-x)}}{e^{-x}-e^{-y}}\\
=& \frac{e^{-(x+\kappa(y-x))}}{e^{-x}}
\frac{1-e^{-(1-2\kappa)(y-x)}}{1-e^{-(y-x)}}
\\
=& e^{-\kappa(y-x)}
\frac{1-e^{-(1-2\kappa)(y-x)}}{1-e^{-(y-x)}}.
\end{align*}
Define $\iota:=\inf_{w\in [0,2c]}e^{-\kappa w}
\frac{1-e^{-(1-2\kappa) w}}{1-e^{-w}}$, then, it can be easily seen that $\iota>0$, and  
$\frac{\mu(L,U)}{\mu(x,y)}>\iota$. Therefore, $A=1-2\kappa$, $B=1$ and $C=\min\{e^{-c}(1-e^{-1}), \iota\}$.

{\hfill$\blacksquare$}\end{exm}

\begin{exm}
\label{exm:exponential_family}
Exponential family distributions with parameter $\ta$ consist  probability measures whose density functions can be represented as $f(x|\ta)= h(x) g(\ta)\exp[\eta(\ta) \tau (x) ]$
with given proper functions $\tau (x)$, $h(x)> 0$, $\eta(\ta)$ and $g(\ta)> 0$, for properties (such as continuity) and examples on these functions, see e.g.~\cite{brown1986fundamentals}.
We assume that the density function is monotone decreasing outside a compact set of $z$, and there exists $M>0$, such that for any $y-x>M(\ta)$, we have,
\begin{align*}
\frac{\int_{x+(1-\kappa) (y-x)}^y h(z) g(\ta)\exp[\eta(\ta) \tau (z)]dz}{\int_x^y h(z) g(\ta)\exp[\eta(\ta) \tau (z)]dz} <\frac{1}{4}.
\end{align*}
Note that this condition holds when the distribution is light-tailed. We can see that 
satisfy the $(A,B,C)$-condition with $A=1-2\kappa$ and $C=e^{-c}$ for some $c>0$ and $\kappa\in (0,\frac12)$; 
\begin{align*}
\frac{\mu(L,U)}{\mu(x,y)} =& \frac{\int_{x+c}^{x+(1-\kappa) (y-x)} h(z) g(\ta)\exp[\eta(\ta) \tau (z)dz}{\int_x^y h(z) g(\ta)\exp[\eta(\ta) \tau (z)dz}
\end{align*}
again can be lower bounded by a constant for properly selected $c$ and $\kappa$.
In the case of $y-x > M(\ta)$, 
\begin{align*}
\frac{\mu(L,U)}{\mu(x,y)} =& \frac{\int_{x+c}^{x+(1-\kappa) (y-x)} h(z) g(\ta)\exp[\eta(\ta) \tau (z)]dz}{\int_x^y h(z) g(\ta)\exp[\eta(\ta) \tau (z)]dz}\\ = 1 & -\underbrace{\frac{\int_x^{x+c}  h(z) g(\ta)\exp[\eta(\ta) \tau (z)]dz}{\int_x^y h(z) g(\ta)\exp[\eta(\ta) \tau (z)]dz}}_{I}
\\
& -\underbrace{\frac{\int_{x+(1-\kappa) (y-x)}^y h(z) g(\ta)\exp[\eta(\ta) \tau (z)]dz}{\int_x^y h(z) g(\ta)\exp[\eta(\ta) \tau (z)]dz}}_{II}.
\end{align*}
Under our assumptions, we can see that both $I$ and $II$ can be uniformly upper bounded by $1/4$. Hence, we can have, $\frac{\mu(L,U)}{\mu(x,y)}\ge \frac12$.
In the case of $y-x\le M(\ta)$, we have, $L(x,y) = x+\kappa(y-x)$ and $U(x,y) = y-\kappa(y-x)$
\begin{align*}
\frac{\mu(L,U)}{\mu(x,y)} =& \frac{\int_{x+\kappa(y-x)}^{y-\kappa(y-x)} h(z) g(\ta)\exp[\eta(\ta) \tau (z)]dz}{\int_x^y h(z) g(\ta)\exp[\eta(\ta) \tau (z)]dz}.
\end{align*}
Under continuity assumption on the density, 
we can see that this is bounded from below by a positive number. 
{\hfill$\blacksquare$}\end{exm}
\begin{lemma}
\label{lem:wd_rec}
If the reference measure $\mu$ satisfies $(A,B,C)$-condition for some real numbers $A<1$ and $B,C >0$,  
then, for any measure $\rho$ on $\Real^+$ and any $p\in (1,\infty)$, we have,  
\begin{align*}
W_p^p (\opTt_\mu\pi, \opTt_\mu\rho)\le W_p^p (\pi, \rho)-C(1-A^p) W_p^{p+B} (\pi, \rho),
\end{align*}
with $\pi$ being the invariant measure.
\end{lemma}
\begin{proof}
Let $\gamma_p(x,y)$ represents the minimum joint distribution for the $W_p$ distance between $\pi$ and $\rho$, i.e.  
\begin{align*}
\int_0^\infty \int_0^\infty |x-y|^p \gamma_p(x,y) dx dy = \inf_{\gamma\in \Gamma(\pi, \rho)} \int_0^\infty \int_0^\infty |x-y|^p \gamma(x,y) dx dy,
\end{align*}
with $\Gamma(\pi, \rho)$ denotes the set of joint distributions on $\RRealp$ with $\pi$ and $\rho$ as the marginals. 
Then,
following the definition in~\eqref{eqndef:gammahat},
for any   $\hat{\gamma}\in\Gamma(\opTs_\mu(\pi), \opTs(\rho))$, we have, 
\begin{align*}
W_p^p& (\opTs_\mu\pi, \opTs_\mu\rho)\le 
\int_0^\infty \int_0^\infty |x-y|^p {\hat \gamma}(x,y) dx dy = 
\\=& 
\int_0^\infty \int_0^\infty \int_0^\infty \mu(z)|x-y|^p {\hat \gamma}(x,y) dx dy\, dz
\\=& 
\int_0^\infty \left(\iint_{I_1}+ \iint_{I_2}+\iint_{I_3}+\iint_{I_4}\right)\mu(z)|x-y|^p {\hat \gamma}(x,y) dx dy\, dz
\end{align*}
with $I_1:=\{x\le z, y\le z\}$, $I_2=\{x\le z, y\ge z\}$, $I_2=\{x\ge z, y\le z\}$ and $I_4=\{x\ge z, y\ge z\}$. Then,
\commentOUT{changed $\opTt$ to $\opTs$, but with index $\mu$ it can be written as simply as $\opT_\mu$}
\commentOUT{Changed $\int_I$ to $\iint_I$}
\begin{align*}
\int_0^\infty \iint_{I_1}&\mu(z)|x-y|^p {\hat \gamma}(x,y) dx dy dz=
\\=
& \int_0^\infty \int_0^z\int_0^z\mu(z)|x-y|^p {\hat \gamma}(x,y) dx dy dz
\\=&
\int_0^\infty \int_0^z\int_0^z\mu(z)|x-y|^p \gamma(z-x,z-y) dx dy dz
\\=&
\int_0^\infty \int_0^z\int_0^z\mu(z)|x-y|^p \gamma(x,y) dx dy dz
\\=& 
\int_0^\infty \int_0^\infty \mu(x\vee y, \infty) |x-y|^p\gamma(x,y) dx dy.
\end{align*}
Similarly 
\begin{align*}
\int_0^\infty \iint_{I_2}&\mu(z)|x-y|^p {\hat \gamma}(x,y) dx dy dz =
\\=&
\int_0^\infty \int_0^\infty\int_0^\infty\mu(z)|x-y|^p  \gamma(z-x,y+z) dx dy dz
\\=& 
\int_0^\infty \int_0^z\int_0^\infty\mu(z)|x-y|^p  \gamma(z-x,y+z) dx dy dz
\\=&
\int_0^\infty \int_0^z\int_z^\infty\mu(z)|2z-x-y|^p  \gamma(x,y) dx dy 
\\=& 
\int_0^\infty \int_0^\infty\left[\int_x^y\mu(z)|2z-x-y|^p  dz\right] \gamma(x,y) dx dy.
\end{align*}
And then 
\begin{align*}
\int\limits_x^y&\mu(z)|2z-x-y|^p  dz =
\\=&  
\left(\int\limits_x^{\al(x,y)}+\int\limits_{\al(x,y)}^{\beta(x,y)}
+ \int\limits_{\beta(x,y)}^y\right) \mu(z)|2z-x-y|^p  dz
\\ 
\le & 
\int\limits_x^{\al(x,y)} |x-y|^p \mu(z) dz + A^p\int\limits_{\al(x,y)}^{\beta(x,y)}|x-y|^p \mu(z) dz + \int\limits_{\beta(x,y)}^y|x-y|^p \mu(z) dz
\\ =& 
|x-y|^p\left(\int_x^y   \mu(z) dz -[1-A^p]\int\limits_{\al(x,y)}^{\beta(x,y)}\mu(z) dz\right).  
\end{align*}
Similar estimation can be obtained for $I_3$ as well.
\begin{align*}
\int_0^\infty \iint_{I_4}&\mu(z)|x-y|^p {\hat \gamma}(x,y) dx dy dz = 
\\ = & \int_0^\infty \int_z^\infty\int_z^\infty\mu(z)|x-y|^p {\hat \gamma}(x,y) dx dy dz 
\\
= &\int_0^\infty \int_0^\infty\int_0^\infty\mu(z)|x-y|^p \gamma(x,y) dx dy dz
\\
= & \int_0^\infty\int_0^\infty |x-y|^p \gamma(x,y) dx dy  
\end{align*}
In summary, we have, 
\begin{align*}
\int_0^\infty \int_0^\infty \int_0^\infty & \mu(z)|x-y|^p {\hat \gamma}(x,y) dx dy dz
\\
 \le &\int_0^\infty \int_0^\infty |x-y|^p \gamma(x,y) dx dy- \\ &- (1-A^p)\int_0^\infty \int_0^\infty |x-y|^p \gamma(x,y) dx dy\int\limits_{\al(x,y)}^{\beta(x,y)}\mu(z) dz.  
\end{align*}
Furthermore, we have, 
\begin{align*}
\int_0^\infty \int_0^\infty &|x-y|^p \gamma(x,y) dx dy\int\limits_{\al(x,y)}^{\beta(x,y)}\mu(z) dz \ge
\\
\ge & 
C  \int_0^\infty \int_0^\infty |x-y|^{p+B} \gamma(x,y) dx dy\\ \ge& C\left(\int_0^\infty \int_0^\infty |x-y|^p \gamma(x,y) dx dy\right)^{\frac{p+B}{p}},
\end{align*}
where the second inequality is an application of the H\"older's inequality. For Wasserstein metric $W_p$, we have, 
\begin{align*}
    W_p^p(\opTs\rho,\pi)\le&W_p^p(\hat{\gamma})\le 
    W_p^p(\gamma)-C(1-A^p)\left(\iint |y-x|^p\gamma(x,y)\,dx\,dy\right)^{1+\frac{B}{p}}\\=&W_p^p(\rho, \pi)-C(1-A^p) W_p^p(\rho,\pi)^{1+\frac{B}{p}}\,.
\end{align*}
{\vskip-0.25cm\hfill$\Box$}\end{proof}

\begin{thm}\label{thm:Wasserstein poly conv}
If for some $A <1 $ and $B,C>0$ the reference measure $\mu$ satisfies $(A,B,C)$-condition then $\opTs^n$ converges at least at a polynomial rate. 

{\hfill$\Box$}
\end{thm}

\section{Self-serving invariant probability}
\label{sec:inv_prob}

In this section, inspired by the Examples~\ref{exm:expon} and~\ref{exm:Bernoulli}, we investigate possible limits of the the iterations of $\opTs_\mu$. In particular, we discuss the intriguing problem of identifying the probability measure that invariant under the transformation $\opT$ with itself also being the reference measure (the fixed point of the transformation defined as $\mu\mapsto \opTs_\mu(\mu)$). In other words, we are trying to identify the random variables $X$ satisfy $|X-X|\stackrel{d}{\sim} X$ (have same distribution). 

\subsection*{Measures with a discrete support}
For a distribution on discrete increasing sequence $x_i$ with  $p_i=p(x_i)_{i=0}^\infty$ the operator $\opT_p(p)$ has the following form. 
\begin{align}\label{eqn:discrete operation}
\hat{p}_k=\opT_p(p)_k=\sum_{|x_i-x_j|=x_k}p_ip_j
\end{align}

\begin{lemma}[Form of the lattice]\label{lem:form of the lattice}
   When $\hat{p}=p$ the support $(x_n)$ of ${p}$ is contained in some initial segment of the lattice L=$x\Integer$,   $(x_n)=(nx)_{n=0}^N$, $0\le N\le \infty$ for some $0<x\in\Realp$. In other words, $x_i=i$ for $i\in[0\dots N]$.
\end{lemma}
\begin{proof}   
By Section~\ref{sec:Trajectories}, discrete support is contained in the lattice $L=z+w\Integer$, for some real $0\le z<w$.
For $i\ge 0$, let $y_i=z+iw$ and $q_i=p(\{y_i\})$. Let $m$ be minimal with $q_m=p(\{y_m\})>0$.  we have
 $\hat{p}(\{0\})\ge q_m^2>0$, thus $0$ is the minimal element of the support of $\hat{p}=p$, $y_m=0=x_0$, $m=0$ and $z=0$ and $q_0=p_0>0$. If the support consists only of $0$, then $p_0=1$ and we are done. If it consists only of two points $x_0=0$ and $x_1=y_k=wk=x_1$ with $p_0+p_1=q_0+q_k=1$,  set $x=wk$ and we are done as well. Otherwise let $0<wi=y_i<y_j=wj$ are in the support of $p$ then $q_{j-i}\ge q_iq_j>0$ and $w(j-i)$ is also in the support. Let $j=di+r$ for some $d,r\in\Integer$. Then for all $m\in\{0,\dots,d\}$ the points $y_{j-mi}$ are in the support with $y_r$ a minimal among them.
Let $x_1=y_k=kw>0$ be the second minimal element of the support. Then for all $y_m$ in the support $k$ divides $m$ the point $0<y_{m\mod k}<y_k$ would be minimal. The point $x=x_1=y_k$ generates the initial segment of the lattice with $N=\sup\{n:y_n\in L\}$.
{\hfill$\Box$}\end{proof}

Therefore for $\hat{p}=p$, Equation~\eqref{eqn:discrete operation} becomes:
\begin{align}
\label{eqn:p_n}
 p_0=\sum_{|i-j|=0}p_ip_j=\sum_{i}p_i^2,\quad  
     {p}_{k}=\sum_{|i-j|=k}p_ip_j=2\sum_{n\ge 0}p_np_{n+k}, k>0.
\end{align}

\begin{exm}\label{exm:half half}[See Example~\ref{exm:Bernoulli}]
The measure $p_0=p_1=\frac{1}{2}$ satisfies $p=\opT_p(p)$ as we have $p_0=p_0^2+p_1^2$ and $p_1=2p_0p_1$.
{\hfill$\blacksquare$}\end{exm}

\begin{exm}\label{exm:geometric self}
    For any $q\in(0,1)$ the distribution $p(q)$ given by $p_0=\frac{1-q}{2}$, $p_k=\frac{1-q^2}{2} q^{k-1}$, for $k>0$ satisfies $\sum_{n\ge 0}p_n=1$ and $\hat{p}(q)=p(q)$.

    Indeed $\sum_{k>0}p_k=\frac{1+q}{2}$ which added to $p_0$ gives 1. Next  summing up geometric series we get
    by Equation~\eqref{eqn:discrete operation}:
    \begin{align*}
    \hat{p}_0&=\sum_{k\ge 0} p_k^2=\left(\frac{1-q}{2}\right)^2+\left(\frac{1-q^2}{2}\right)^2\frac{1}{1-q^2}=\frac{2-2q}{4}=p_0
\\
    \hat{p}_n&=2p_0p_n+2\sum_{k>0}p_kp_{n+k}=
    (1-q)p_n+2\left(\frac{1-q^2}{2}\right)^2q^{n-1}q\sum_{k>0}(q^2)^{k-1}\\
    &=(1-q)p_n+\frac{1-q^2}{2}q^{n-1} q=(1-q)p_n+p_n q=p_n
    \quad\text{so that}\quad\hat{p}=p\,. 
   \end{align*}
{\vskip-0.25cm\hfill$\blacksquare$}\end{exm}

\begin{prop}[$p=\hat{p}$, the discrete case]
If $p_0=\mu(0)=1$ then $\hat{p}=p$. 
Assume that $p=\hat{p}$ and $p$ has a discrete support.
\begin{enumerate}
\item 
    When $p_0\not=1$ then $0<p_0\le \frac{1}{2}$. 
\item 
    When $p_0\not=1$ and  only a finite number of $p_i\not=0$  then $p_0=\frac{1}{2}$. 
\item 
    If $p_0=\frac{1}{2}$ then $p_1=\frac{1}{2}$ and all the other $p_i=0$, $i>1$.        
\end{enumerate}
\end{prop}

\begin{proof}
  If $p_0=1$ and all the other  $p_i=0=\hat{p}_i$ 
  then by Equation~\eqref{eqn:discrete operation} $\hat{p}_0=p_0^2=1$ and all the other  $\hat{p}_i=0=p_i$ which proves invariance.
By Lemma~\ref{lem:form of the lattice},  we have $p_0>0$.
From now on we assume  $\hat{p}=p$ and $p_0<1$.
\begin{enumerate}
 \item 
If $p_0<1$ then $p_k>0$ for some $k>0$, and $0<p_k=p_0p_k+p_kp_0+2\sum_{n>0}p_np_{n+k}\ge p_k(2p_0)$, that is  $1\ge 2p_0$.
\item
  Suppose $0<k<\infty$ is largest index  with $p_k>0$. Then $p_k=\hat{p}_k=p_0p_k+p_kp_0$, as these are the only pairs of indices to provide $|i-j|=k$ with positive probability. As $p_k>0$ we get  $1=2p_0$. 
  \item
    If $p_0=\frac{1}{2}$ then for $k\not=0$ $p_k=\hat{p}_k=p_0p_k+p_kp_0+\sum_{n\not=0} (p_{n+k} p_n+p_np_{n+k})=
    p_k+2\sum_{n>0}p_np_{k+n}$. Then the sum $2\sum_{n>0}p_n p_{k+n}=0$ for every $n,k>0$. Suppose there are two distinct $i,j>0$ with $p_i,p_j>0$ then $p_ip_j>0$ and the sum would be positive. That proves uniqueness of $k>0$ with $p_k>0$ and it must satisfy $1=p_0+p_k=\frac{1}{2}+p_k$. All the other $p_i=0$ and by the ordering of the support $x_i$, we have $k=1$.
\end{enumerate}
{\hfill$\Box$}\end{proof}

\begin{coro}\label{cor:types of support}
The support of $p$ is either   $\{0\}$, $\{0,x\}$ or the whole $x\Natural$.
\end{coro}

\subsection*{Measures with a density in $L^1(dt)$}
Let's assume that the reference measure $\mu$ has a density with respect to the Lebesgue measure so that (abusing the notation again) $\mu(dt)=\mu(t)dt$ with $\mu(t)\in L^1(dt)$. There are no atoms of $\mu$, in particular there is no atom at 0 and for any Lebesgue measurable function $f:\Realp\times\Realp\to\Real$ the double integral over a set $A\subset\Realp\times\Realp$ of zero Lebesque measure $\iint_A dtds=0$, yields $\iint_A f(t,s)\mu(dt)\mu(ds)=0$.
Suppose that  $\mu=\opTs_\mu(\mu)$.
Then
\begin{align*}
    \mu(x)&=\int_{t\ge 0}\mu(t+x)\mu(t)\,dt+\int_{t\ge 0}\mu(t)\mu(t+x)\,dt=2\int_{t>0}\mu(t+x)\mu(t)\,dt\,.
\end{align*}
\begin{exm}\label{exm: continuous example}[See Example~\ref{exm:expon}]
    Let $\mu(t)=\alpha e^{-\alpha t}$, for $t\ge 0$ and $\alpha>0$. 
    Then 
    \begin{align*}
        \hat{\mu}(x)&=2\int_0^\infty\alpha e^{-\alpha t}\alpha e^{-\alpha(t+x)}\,dt
        =
        \alpha e^{-\alpha x}\int_0^\infty2\alpha e^{-2\alpha t}\,dt=\alpha e^{-\alpha x}=\mu(x)
    \end{align*}
{\vskip-0.25cm\hfill$\blacksquare$}\end{exm}
\subsection*{The form of generating functions for self-serving distributions}
\subsubsection*{The discrete case}
Let for $z\in\Complex$, $f(z)=\sum_{n=0}^\infty p_n z^n$ be a probability generating function of $\mu$, and $g(z)=2f(z)-1$.
\begin{exm}\label{exm:geometric self f g}
    In  Example~\ref{exm:geometric self}, we have
    $f(z)=\frac{1-q}{2}\left(\frac{1+z}{1-qz}\right)$ and $g(z)=\frac{z-q}{1-qz}$, which follows from an elementary calculation on geometric series.
{\hfill$\blacksquare$}\end{exm}

\begin{prop}
Both functions $f$ and $g$ are analytical in the unit disk $D^o=\{x:|z|<1\}$ and continuous in its closure $\overline{D}=\{z:|z|\le 1\}$.
On the boundary, \emph{i.e.} on the unit circle $|z|=1$, we have
\begin{align*}
f(z) \overline{f(z)}=\frac{1}{2} \left(f(z) + \overline{f(z)}\right),\quad  
g(z)\overline{g(z)}=1,\quad 
|g(z)|=1\,.
\end{align*}
\end{prop}
\begin{proof}
The case $p_0=1$ and $p_0=p_1=\frac{1}{2}$ was presented in Example~\ref{exm:half half}, here $g(z)=z$. By Corollary~\ref{cor:types of support} it remains to deal with the case when $p_n>0$ for all $n$. 

For $|z|\in \overline{D}$:  we have $|f(z)|\le \sum_{n=0}^\infty p_n |z^n|\le 1$.
The power series for $f'(z)=\sum_{n=1}^\infty n p_n z^{n-1}$ converges absolutely for $z\in D^o$.
That makes both $f$ and $g$ holomorphic in the open unit disk and continuous on the closed unit disk.   then $|z|=1$, Equation~\eqref{eqn:p_n} together with $\hat{p}=p$ implies
\begin{align*}
f(z) \overline{f(z)}&=\sum_{i\ge 0}p_i z^i\sum_{j\ge 0} p_j \overline{z}^j=\sum_{n\ge 0}\sum_{k=0}^n p_k p_{n-k}z^k\overline{z}^{n-k}\\
&=
\sum_{n\ge 0} p_n^2 z^n\overline{z}^n+
\sum_{n\ge 0}\sum_{k>0}p_n p_{n+k}(z^{n+k}\overline{z}^{n}+
z^{n}\overline{z}^{n+k})\\
&=
\sum_{n\ge 0} p_n^2 |z|^{2n}+
\sum_{n\ge 0}|z|^{2n}\sum_{k>0}p_n p_{n+k}(\overline{z}^{k}+
z^{k})\\
&=\sum_{n\ge 0} p_n^2+\sum_{k>0}\left(\overline{z}^{k} + z^{k}\right)\sum_{n\ge 0}p_n p_{n+k}\\
&=
\frac{1}{2} \left(\hat{p}_0+\sum_{k>0}z^k\hat{p}_k +\hat{p}_0+\sum_{k>0}\overline{z}^k\hat{p}_k\right) 
=\frac{1}{2} \left(f(z) + \overline{f(z)}\right),
\end{align*}
where for the absolutely convergent series we used the equalities
\begin{align*}
\sum_{i\ge 0} a_i\sum_{j\ge 0}b_j=&\sum_{n\ge 0}\sum_{0\le k\le n}a_k b_{n-k}\\=&\sum_{i=j}a_ib_j+\sum_{i\ge 0,j>i}a_ib_j+\sum_{j\ge 0, i>j}a_ib_j\\=&\sum_{n\ge 0}a_nb_n+\sum_{n\ge 0,k>0}a_nb_{n+k}+\sum_{n\ge 0, k>0}a_{n+k}b_n.
\end{align*}
Finally, 
$g(z)\overline{g(z)}=4f(z)\overline{f(z)}-2\left(f(z)+\overline{f(z)}\right)+1=1$, and as $|g(z)|=|\overline{g(z)}|=|g(\overline{z})|$, we get  $|g(z)|^2=1$ on the unit circle.
{\hfill$\Box$}\end{proof}

\begin{coro}
The function $g(z)$ has real coefficients and it is a finite Blaschke product in form
\begin{align*}
g(z)&=e^{\ta i} z^m \prod_{n=1}^N \frac{|a_n|}{a_n}\frac{z-a_n}{1-z{\bar a}_n},
\end{align*}
for some real $\ta$, integers $m,N\ge 0$, and complex $\{a_1, a_2, \ldots, a_N\}$. 
\end{coro}
\begin{proof}
As $g$ is holomorphic and $|g(z)|=1$ on the unit circle then by the maximum principle $|g(z)|< 1$ on $D^o$. So $g:\overline{D}\to \overline{D}$ is continuous. 
Therefore, according to a result attributed to Fatou~\cite{FatouBlaschke,garcia2018finite}, it is a finite Blaschke product. 
{\hfill$\Box$}\end{proof}

\begin{rem}
Fatou's argument~\cite{FatouBlaschke,garcia2018finite} is quite simple, consider the two ratios between $g$ and the finite Blaschke product formed by all the zeros of $g$, they are both analytic in $D$ and their norms are one at the boundary. By maximum principle, both ratios are bounded by one. Hence, the ratios have to be a constant one.   
\end{rem}
\subsubsection*{The continuous case}
For $z\in\Complex$
(with $r,\zeta \in\Real$) define
\[
{z} \mapsto f(z)=f_\mu(z)=\int_{t\ge 0}e^{{z} t}\mu(t)\,dt\,.
\] 
The function $f$ as a $\mu$ weighted average of holomorphic functions $z\to e^{-zt}$, $t\ge 0$, is holomorphic itself at least on the half-plane $\Re(z)<0$ and $f(0)=1$.
As before $g(z)=2f(z)-1$ is a holomorphic function in the same domain as $f$. 
\begin{exm}\label{exm: continuous example f g}
In Example~\ref{exm: continuous example} we have
$f(z)=\frac{\alpha}{\alpha-z}$ and $g(z)=\frac{\alpha+z}{\alpha-z}$, which follows by (careful but) elementary calculation.
{\hfill$\blacksquare$}\end{exm}
\begin{align*}
    g(z)=\frac{z-\alpha+2\alpha}{\alpha-z}=-1+2\frac{1}{1-\frac{z}{\alpha}}=-1+2\sum_{n\ge 0}\left(\frac{z}{\alpha}\right)^n=1+2\sum_{n>0}\left(\frac{z}{\alpha}\right)^n
\end{align*}
\begin{lemma}[Formulae for $f$ and $g$ when $\mu\in L^1(dt)$]
\label{lem:formule for mu in L1}
For 
\begin{align*}
    f(z )\cdot f(-z )&=\frac{1}{2}\left(f(z )+f(-z )\right)\\
    g(z )g(-z )&=1
\end{align*}   
\end{lemma}

\begin{proof}
As the Lebesgue measure of the sets $\{(t,s):t=s\}$, $\{(0,s)\}$ and  $\{(t,0)\}$  is zero  0 in $\Realp\times\Realp$ we have
\begin{align*}
&f(-z )\cdot f(z )=\int\limits_{t\ge 0}e^{-{z } t}\mu(t)\,dt \cdot \int\limits_{s\ge 0} e^{{z } s}\mu(s)\,ds
\\
=&\iint\limits_{t\ge 0,s\ge 0}e^{-{z } t}\mu(t)  e^{{z } s}\mu(s)\,dtds
\\
=&\iint\limits_{0<s<t}e^{-{z } (t-s)}\mu(t) \mu(s)\,dtds+\iint\limits_{0<t<s}e^{{z } (s-t)}\mu(t) \mu(s)\,dtds\\
=&
\iint\limits_{s> 0, w> 0}e^{-{z } w}\mu(w+s)\mu(s)\,dwds
+
\iint\limits_{t> 0, w> 0}e^{{z } w}\mu(t)\mu(w+t)\,dtdw\\
=&
\int\limits_{w>0}e^{-{z } w}(\int\limits_{s>0}\mu(w+s)\mu(s)\,ds)\,\,dw+
\int\limits_{w>0}e^{{z } w}(\int\limits_{t>0}\mu(t)\mu(w+t)\,dt)\,\,dw 
\\
=&
\frac{1}{2}\int\limits_{w>0}e^{-{z } w}\mu(w)\,dw+\frac{1}{2}\int\limits_{w>0}e^{{z } w}\mu(w)\,dw 
=\frac{f(z )+f(-z )}{2}\,.
\end{align*} 
The proof  of $g(z)\cdot g(-z)=1$ is similar as in the discrete case.
{\hfill$\Box$}\end{proof} 

\begin{coro}\label{cor:g on boundary}
    On the boundary $\Re(z)=0$ we have $-z=\bar{z}$ and there 
\begin{align*}
    f(z )\cdot f(\bar{z} )&=\frac{1}{2}\left(f(z )+f(\bar{z})\right)\\
    g(z )g(\bar{z})&=1
\end{align*}   
and as $|g(z)|=|g(\bar{z}|$ we get $|g(z)|=1$ for $\Re(z)=0$.
\end{coro}
Let $\zeta:\Complex\to\Complex$, $z\mapsto\zeta(z)=\frac{z-1}{z+1}$. The unit circle is mapped onto the imaginary axes and its interior to the left half plane $\Re(z)<0$. In particular  $0\mapsto -1$, $1\mapsto 0$, $\pm i\mapsto\pm i$.
\begin{rem}
    One can check that $f$ and $g$ from Example~\ref{exm:geometric self f g} are the same as $f\circ\zeta$ and $g\circ\zeta$ from Examples~\ref{exm: continuous example f g} when for $0<\alpha<\infty$ and $1>q>-1$ one uses the relation $\alpha q+\alpha+ q=1$.
\end{rem}
\begin{prop}
The holomorphic map $h=g\circ \zeta$ maps the unit disk to the  unit disk, the unit circle to the unit circle and therefore it is a finite Blaschke product. 
\end{prop}
\begin{proof}
    Combine Lemma~\ref{lem:formule for mu in L1} and Corollary~\ref{cor:g on boundary}.
{\hfill$\Box$}\end{proof}
\subsection*{Open problems}
    We do not know if there exist a Blaschke product with more than one factor which correspond to a probability $\mu$ satisfying $\opTs_\mu\mu=\mu$, in other words if our examples are all possible such probabilities. 

    We also do not know the behavior of the  iterates of the non-linear operator on the space of probability measures $\mu\mapsto \opTs_\mu\mu$. Is there only a one parameter family of fixed points? Are they attracting? If there is a convergence,  what feature of the starting measure laminates the space into the bassins of different atractors?

\bibliographystyle{plain}

\end{document}